\documentclass[preprint,11pt]{elsarticle}

\usepackage{lineno}
\usepackage[hidelinks]{hyperref}

\usepackage{kotex}
\usepackage{acronym}
\usepackage{mathtools}
\usepackage{algorithmic}
\usepackage{algorithm}
\usepackage{amssymb}
\usepackage{nccmath}

\usepackage{amsthm}
\usepackage{dsfont}
\usepackage{dirtytalk}
\usepackage{afterpage}
\usepackage{ulem}
\usepackage{bm}
\usepackage{longtable}
\usepackage{caption} 
\usepackage{hhline}
\usepackage{lscape}
\usepackage{array,multirow}
\usepackage{makecell}
\usepackage{array}
\usepackage{xcolor}
\usepackage{circledsteps}
\usepackage{tikz}

\makeatletter
\newtheorem*{rep@theorem}{\rep@title}
\newcommand{\newreptheorem}[2]{%
\newenvironment{rep#1}[1]{%
 \def\rep@title{#2 \ref{##1}}%
 \begin{rep@theorem}}%
 {\end{rep@theorem}}}
\makeatother

\newtheorem{proposition}{Proposition}
\newreptheorem{proposition}{Proposition}
\newtheorem{definition}{Definition}

\usepackage{setspace}

\usepackage{graphicx}
\usepackage{subcaption}
\usepackage{mwe}

\journal{European Journal of Operational Research}

\biboptions{authoryear}

\begin{document}

\begin{frontmatter}

\title{Effectiveness of Social Distancing under Partial Compliance of Individuals}



\author[mymainaddress]{Hyelim Shin}
\ead{hyelim.shin@kaist.ac.kr}

\author[mymainaddress]{Taesik Lee}
\ead{taesik.lee@kaist.edu}


\address[mymainaddress]{Department of Industrial and Systems Engineering, KAIST, 291 Daehak-ro, Yuseong-gu, Daejeon, Republic of Korea}

\begin{abstract}

Social distancing reduces infectious disease transmission by limiting contact frequency and proximity within a community. However, compliance varies due to its impact on daily life. This paper explores the effects of compliance on social distancing effectiveness through a "social distancing game," where community members make decisions based on personal utility. We conducted numerical experiments to evaluate how different policy settings for social distancing affect disease transmission.

Our findings suggest several key points for developing effective social distancing policies. Firstly, while generally effective, overly strict policies may lead to noncompliance and reduced effectiveness. Secondly, the public health benefits of social distancing need to be balanced against social costs, emphasizing policy efficiency. Lastly, for diseases with low reinfection risk, a segmented policy exempting immune individuals could lessen both infections and socioeconomic costs.
\end{abstract}

\begin{keyword}
OR in health services \sep epidemics \sep social distancing \sep compliance \sep game theory
\end{keyword}

\end{frontmatter}

\section{Introduction}
\label{sec:INTRO}

An epidemic outbreak can severely strain healthcare systems, sometimes leading to a complete breakdown, as observed during the COVID-19 crisis. The primary objective for public health authorities in such situations is to contain the disease spread so that the existing healthcare system can handle the demand from the outbreak while still maintaining its function for regular services. While vaccination is the most effective defense against an epidemic outbreak, it's often unavailable in the early stages due to the extensive time and resources required for development, approval, and distribution.

As an alternative, non-pharmaceutical interventions (NPIs) including social distancing are often employed. Social distancing reduces disease transmission by limiting frequency and proximity of contacts between people \citep{CDC-socialdistancingdef}. While its success largely depends on community compliance, there is a varying degree of compliance among community members as it can significantly disrupt daily life \citep{nivette2021non-noncompliance,hills2021factors-noncompliance,rosha2021factors-noncompliance,yamamoto2021quantifying-noncompliance,hengartner2022factors-noncompliance}. For instance,  \cite{hills2021factors-noncompliance} show that 92.8$\%$ of participants in their study admitted to violating distancing rules during COVID-19. Therefore, understanding and incorporating community compliance is crucial when designing and implementing these policies.

This study explores how compliance with social distancing influences its effectiveness in controlling disease spread, by introducing a ``social distancing game" into an epidemic model. In this game, community members decide on their compliance with the social distancing policy based on personal utility. This decision problem is cast as a game model because their personal utility is affected by other members' compliance decision. Our numerical experiments show that for social distancing to be effective, the policy design must take into account the effect of possible noncompliance and include measures to encourage compliance with the policy. \textcolor{black}{
Specifically, our work makes three primary contributions. First, we introduce a novel epidemic model that incorporates a social distancing game, integrating individual responses with social distancing policies. Second, we identify conditions for compliance by calculating the Nash equilibrium for the game, revealing that overly stringent policies can lead to noncompliance and hence necessitate substantial penalties for effective enforcement. Third, our experimental results offer suggestions for more effective social distancing policies. For example, we evaluate the social and economic impacts of strict social distancing measures and the potential benefits of a segmented social distancing strategy.}

The structure of this paper is as follows: Section \ref{sec:literature} reviews epidemic models incorporating social distancing. Section \ref{sec:problem} outlines the social distancing game, while Section \ref{sec:solution} examines compliance decisions via Nash equilibrium and Section \ref{sec:result} presents experimental results under different policy scenarios. Additional insights into enhancing policy efficiency are discussed in Section \ref{sec:disc}. The paper concludes with a summary of findings in Section \ref{sec:conclusion}.


\section{Literature Review}
\label{sec:literature}

Our paper introduces a novel mathematical model that captures the dynamics of infectious disease spreading under social distancing policies. Mathematical epidemic models, often referred to as a compartmental model, are originated from pioneering work by \cite{kermack1927-SIR}. These models use compartments to represent different population groups -- such as susceptible, infected, and recovered -- and transition rates to depict the dynamics of disease spread. They provide a simplified framework for analyzing the dynamics of disease spread on a large scale, making compartmental models ideal for modeling national social distancing policies.

Several studies have extended the basic compartmental model to account for social distancing effects. For example, \cite{brauer2011-BSIR} and \cite{zhong2013-imit} adjust the transmission rate to reflect social distancing measures, with populations divided into avoidance and non-avoidance groups. Other models, such as those by \cite{liu2015-func} and \cite{nyabadza2010-func}, introduce functions that adjust individual behaviors in response to the infected population's prevalence, positing that a higher infection rate leads to increased social distancing.

Additionally, some studies incorporate external influences on social distancing behavior. \cite{poletti2012-imit} explore how imitation impacts behavior, while \cite{eryarsoy2023models-func} consider the cost of interventions in their transmission rate adjustments. The influence of media is also significant, with studies like those by \cite{perra2011-func} and \cite{greenhalgh2015-media} examining how media campaigns and public awareness affect social distancing practices.

An alternative modeling approach treats social distancing as a utility-maximizing decision within a game-theoretical framework. \cite{fenichel2011-BSIR} and \cite{reluga2010-BSIR} model daily contact decisions and vaccination costs, respectively, as strategic games among individuals. Similarly, \cite{farboodi2021-BSIR} and \cite{quaas2021-BSIR} use optimal control models to analyze behaviors under different policy scenarios, such as laissez-faire and altruism.

While previous studies have primarily examined how social distancing behavior changes in response to the speed or magnitude of transmission of infectious diseases or the behavior of those around them, few have investigated the effects from a social distancing policy itself. To address this gap, we propose a novel social distancing game model that captures the impact of people's decisions in response to a specific social distancing policy. In particular, our model considers a social distancing policy that includes both public guidance for the intensity of social distancing and penalties for noncompliance, and examine how a transmission rate changes by individuals' behavioral decisions.

To the best of our knowledge, \cite{aurell2022optimal-BSIR(wPolicy)} is the only prior study that models a feedback from a social distancing policy on the individuals' social distancing behavior. They define a social distancing policy in terms of a regulator's decision on the social distancing level and terminal utility, with players who conform more closely to the regulator's socialization level incurring lower costs. However, their model does not consider the utility individuals obtain from social activities. This utility from social activities is the reason why individuals would not simply reduce their social activity level, and our model takes into account the behaviors of individuals to earn utilities from engaging in social activities.


\section{Problem Definition and Formulation}
\label{sec:problem}

\subsection{Problem Description}
\label{subsec:problem-description}

We consider the spread of a novel infectious disease within a closed community. We focus on a period when no vaccine is available, leaving public health authorities to rely solely on social distancing as a mitigation strategy. In our study, we define a social distancing policy using two parameters \((\tilde{\alpha}_t, c_t)\). The parameter \(\tilde{\alpha}_t\) denotes the maximum allowable social activity level; normal social activity is indicated by \(\tilde{\alpha}_t = 1\), and \(\tilde{\alpha}_t = 0\) represents full lockdown policy. A social distancing policy comes with a means for enforcement, which is represented by the penalty for policy violations, \(c_t\). 

We use a compartment model to describe disease spreading. In this model, individuals in the community belong to one of \textit{susceptible, infected, quarantined,} and \textit{recovered} compartments. We assume that individuals in the infected compartment do not yet know they are infected, hence continuing to engage in social contacts, whereas those in the quarantined know of their infectiousness and are disconnected from the rest of the community. 

We investigate the individuals' decision on compliance with a given social distancing policy when they choose their social activity levels to maximize their utility. We define the individual's utility to consist of the pleasure from social activity, risk of infection, and penalties for noncompliance. Individuals' decisions are interdependent since the collective level of social activities determines the degree of disease spreading, thereby one's risk of infection. We represent these interactions through a multi-agent game, which we refer to as a social distancing game. By analyzing this game, we assess how a social distancing policy affects individuals' compliance decisions and how those decisions influence the effectiveness of the social distancing strategy.




\subsection{Epidemic Model with Social Distancing Behavior}
\label{subsec:problem-epidemic}

To account for the compliance behavior in disease spread dynamics, we extend a standard SIR model by including a quarantine compartment. 
It is defined by the following  ordinary differential equations:
\textcolor{black}{
\begin{align}
    \frac{ds(t)}{dt} &= -\beta \frac{s(t)i(t)}{s(t)+i(t)+r(t)} \, ;& \quad
    \frac{di(t)}{dt} &= \beta  \frac{s(t)i(t)}{s(t)+i(t)+r(t)} - \gamma_1 i(t) \, ; \nonumber \\ 
    \frac{dq(t)}{dt} &= \gamma_1 i(t) - \gamma_2 q(t) \, ;& \quad
    \frac{dr(t)}{dt} &= \gamma_2 q(t) \,. \label{eq:SIR}
\end{align}
}

In the model equations, $s(t)$, $i(t)$, $q(t)$, and $r(t)$ represent the susceptible, infected, quarantined, and recovered populations at time $t$, expressed as fractions of the total population $N$. The transmission rate parameter $\beta$ is the product of the nominal contact rate $\kappa$ and the probability of infection per contact $\pi$. \textcolor{black}{Note, in the right-hand-side of the first equation for $ds(t)/dt$, that $q(t)$ is excluded in the denominator as they do not engage in contact with other individuals.} 
The rate at which infected individuals move to the quarantined state is $\gamma_1$, and the recovery rate for quarantined individuals is $\gamma_2$. As typical in epidemic models, these parameters are assumed constant throughout the disease spread.

To incorporate the effect of social distancing, we modify the rate of infection, $-ds(t)/dt$, such that this rate changes according to the social distancing behavior of the population \citep{brauer2011-BSIR, fenichel2011-BSIR}. By rewriting the rate of infection, $-ds(t)/dt$, with $\beta$ separated into $\kappa$ and $\pi$, we have
\begin{align}
    -\frac{ds(t)}{dt} &= \kappa i(t) \times \frac{s(t)}{s(t)+i(t)+r(t)} \times \pi 
    \label{eq:dsdt}
\end{align}
The three terms represent the contacts per unit time by infected individuals, the probability that such a contact is with a susceptible individual, and the probability of transmission upon infectious contact. Then, we modify Eq.\eqref{eq:dsdt} to incorporate social activity levels as follows:
\begin{align}
    -\frac{ds(t)}{dt} 
    &= \, \kappa \alpha^I(t) i(t) \times \frac{\alpha^S(t) s(t) }{\alpha^S(t) s(t) + \alpha^I(t) i(t) + \alpha^R(t) r(t)  } \times \pi \, . \label{eq:adjusted_dsdt}
\end{align}
Time-varying parameters $\alpha^{(\cdot)}(t)$ represent each population group's compliance with social distancing policies. As with $\tilde{\alpha}_t$, these parameters range from 0 to 1, where 1 indicates normal activity level and 0 reflects no social activity at all. 
Then, the first term quantifies the effective number of contacts per unit time by infected individuals, and the second term calculates the probability of these contacts involving a susceptible based on the reduced activity levels of each group as dictated by $\alpha^{(\cdot)}(t)$. 
\textcolor{black}{With this modification, we have an  epidemic model with social distancing behavior:}
\textcolor{black}{
\begin{align}
    \frac{ds(t)}{dt} 
    &= \, -\beta \frac{\alpha^S(t) s(t) \alpha^I(t) i(t)  }{\alpha^S(t) s(t) + \alpha^I(t) i(t) + \alpha^R(t) r(t) } \, ;& \\
    \frac{di(t)}{dt} 
    &= \, \beta \frac{\alpha^S(t) s(t) \alpha^I(t) i(t)  }{\alpha^S(t) s(t) + \alpha^I(t) i(t) + \alpha^R(t) r(t)  } - \gamma_1 i(t) \, ;&  \nonumber \\
    \frac{dq(t)}{dt}
    &= \, \gamma_1 i(t) \, - \gamma_2 q(t);& \nonumber \\
    \frac{dr(t)}{dt} &= \, \gamma_2 q(t) \, . \nonumber 
    \label{eq:adjusted_SIR}
\end{align}
}
\subsection{Objective function for Individuals}
\label{subsec:problem-individual}

Compliance with a social distancing policy is determined as individuals adjust their social activity levels, \(\alpha_t^j\), to maximize utility. To develop the decision-making model, we make the following assumptions. First, individuals are aware of the upcoming epidemic states in the community and social distancing policy requirement: \(\{\hat s(t+1), \hat i(t+1), \hat q(t+1), \hat r(t+1), \tilde{\alpha}_{t+1}, c_{t+1}\}\). This assumption is not unrealistic as social distancing measures are typically announced in advance, and predictions about future epidemic states often guide decision-making. Second, individuals within a compartment exhibit homogeneous behaviors driven by identical utility functions. Though it is possible to further divide each compartment according to utility variations, analyses and interpretations of the results are much more useful by avoiding unnecessary complications. 

As specified in the problem description, we construct a utility function with three main factors: the pleasure from social activities, penalty for noncompliance, and perceived risk of infection. Note that, in our model, compliance decisions by the infected individuals are the same as those by the susceptible since they do not know of their infectiousness and behave as a susceptible. In terms of decision model, this means that these two groups share the same utility function.  \color{black} Thus, the utility maximization objective of individual $j$ in the susceptible and infected group is given as follows:




\begin{equation}
\label{eq:obj-SI}
    \max_{\alpha^j \in [0,1]} (\log\alpha^j-\alpha^j+1)-(c \, \mathds{1}(\alpha^j>\tilde\alpha))-\alpha^j \frac{\alpha^{-j} \, \hat i \, \beta \, d}{\alpha^{-j} (\hat s \,+\, \hat i) \,+\, \alpha^{R} \hat r}.
\end{equation}

\color{black}
\noindent For simplicity, we omit the time index $t$ in $\alpha$ and $c$. This utility function balances the benefits of social interaction against the risks of infection and penalties for noncompliance. The first term of in the objective function measures the utility from social activities, modeled using the increasing concave function from \cite{farboodi2021-BSIR}. The second term represents the potential penalty for noncompliance, which occurs if an individual's social activity exceeds the mandated level, $\tilde{\alpha}$, with $c$ reflecting the penalty amount and detection likelihood. The final term assesses the perceived risk of infection, combining the probability of infection and sensitivity to it, denoted by $d$. \textcolor{black}{To denote the average activity level of all individuals in the susceptible and infected groups, excluding the decision maker $j$, we use $\alpha^{-j}$. In the denominator, we apply $\alpha^{-j}$ to both susceptible and infected group as their decisions are identical. }

\color{black}
Individuals who are recovered from the disease do not have the risk of getting infected again, and thus the third term in Eq.\eqref{eq:obj-SI} becomes irrelevant. By dropping the third term, we obtain the objective function for the recovered population as follows:

\begin{equation}
\label{eq:obj-R}
    \max_{\alpha^j \in [0,1]} (\log\alpha^j-\alpha^j+1)+(-c \mathds{1}(\alpha^j>\tilde\alpha)).
\end{equation}
While we assume no re-infection for the recovered population in Eq.\eqref{eq:obj-R}, later in our numerical experiments, we consider the possibility of re-infection  for the recovered population and examine the effect of social distancing policy in such scenarios.

\section{Individuals' Decision for a Given Social Distancing Policy}
\label{sec:solution}

\subsection{Optimal decision for recovered individuals}
\label{subsec:solution-R}

\begin{figure}[htb]
\centering
\includegraphics[width=0.4\textwidth]{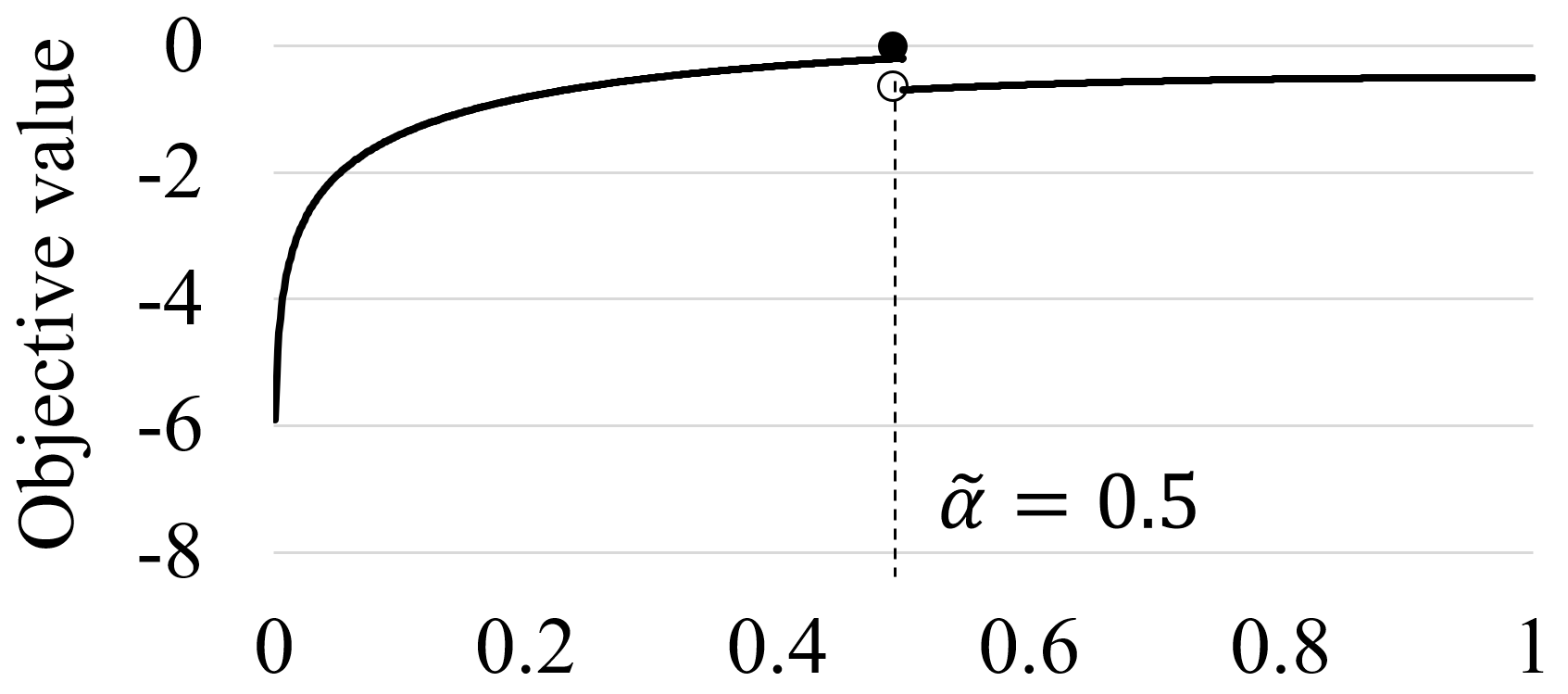}
\caption{Objective function of recovered individuals ($\tilde{\alpha}=0.5$)}
\label{graph:obj-R}
\end{figure}

We first look at the optimal decision for the recovered individuals as they are much simpler than the susceptible and infected. The objective function for recovered individuals in equation Eq.\eqref{eq:obj-R} is plotted in Figure \ref{graph:obj-R}. As shown in the figure, it is monotonically increasing in the intervals $\alpha^j \in [0,\tilde{\alpha}]$ and $\alpha^j \in (\tilde{\alpha},1]$ with a drop at $\alpha^j = \Tilde{\alpha}$. The drop is due to the indicator function in the penalty term: $\mathds{1}(\alpha^j>\tilde\alpha)$. 

This pattern of monotonic increase, with a drop by constant $c$ at $\alpha^j = \tilde{\alpha}$, holds true for any value of $\tilde{\alpha}$. Leveraging this characteristic of the objective function, we can obtain an optimal solution simply by comparing the objective value at $\alpha^j = \tilde{\alpha}$ with the objective value at $\alpha = 1$. Then the optimal solution $\alpha^{*R}$ can be written as
\begin{equation}
\label{eq:decision-R}
    \alpha^{*R} = \begin{cases}
            \tilde{\alpha}, & \text{$\log{\tilde{\alpha}}-\tilde{\alpha}+1>-c$}\\
      1, & \text{otherwise.}
         \end{cases}
\end{equation}

We obtain this solution with ease because the objective function is solely dependent on its own decision $\alpha^R$, and is not influenced by the decisions of other individuals. In contrast, for susceptible and infected individuals, their utility functions are influenced by other individuals' decisions due to the infection risk term. This is discussed in the following section.

\subsection{Optimal decision for susceptible and infected individuals}
\label{subsec:solution-SI}

The last term in the objective function, Eq.\eqref{eq:obj-SI}, for the susceptible and infected  denotes the risk of infection, and it includes $\alpha^{-j}$ which is the social activity level by other individuals in the susceptible and infected groups. 
Consequently, the decision of a susceptible (and infected) individual $j$ on their level of social activity, $\alpha_j$, is influenced by the decisions of others in those groups. This interaction defines the decision-making process for these individuals as a social distancing game. In this game, all participants are assumed rational and act in accordance with the Nash equilibrium concept. Essentially, each individual optimizes their own objective function based on the decisions of others, culminating in a Nash equilibrium where no one can better their situation by changing their strategy alone. 

\begin{figure}[htb]
 \centering
        \begin{subfigure}[b]{0.32\textwidth}
            \centering
            \includegraphics[width=\textwidth]{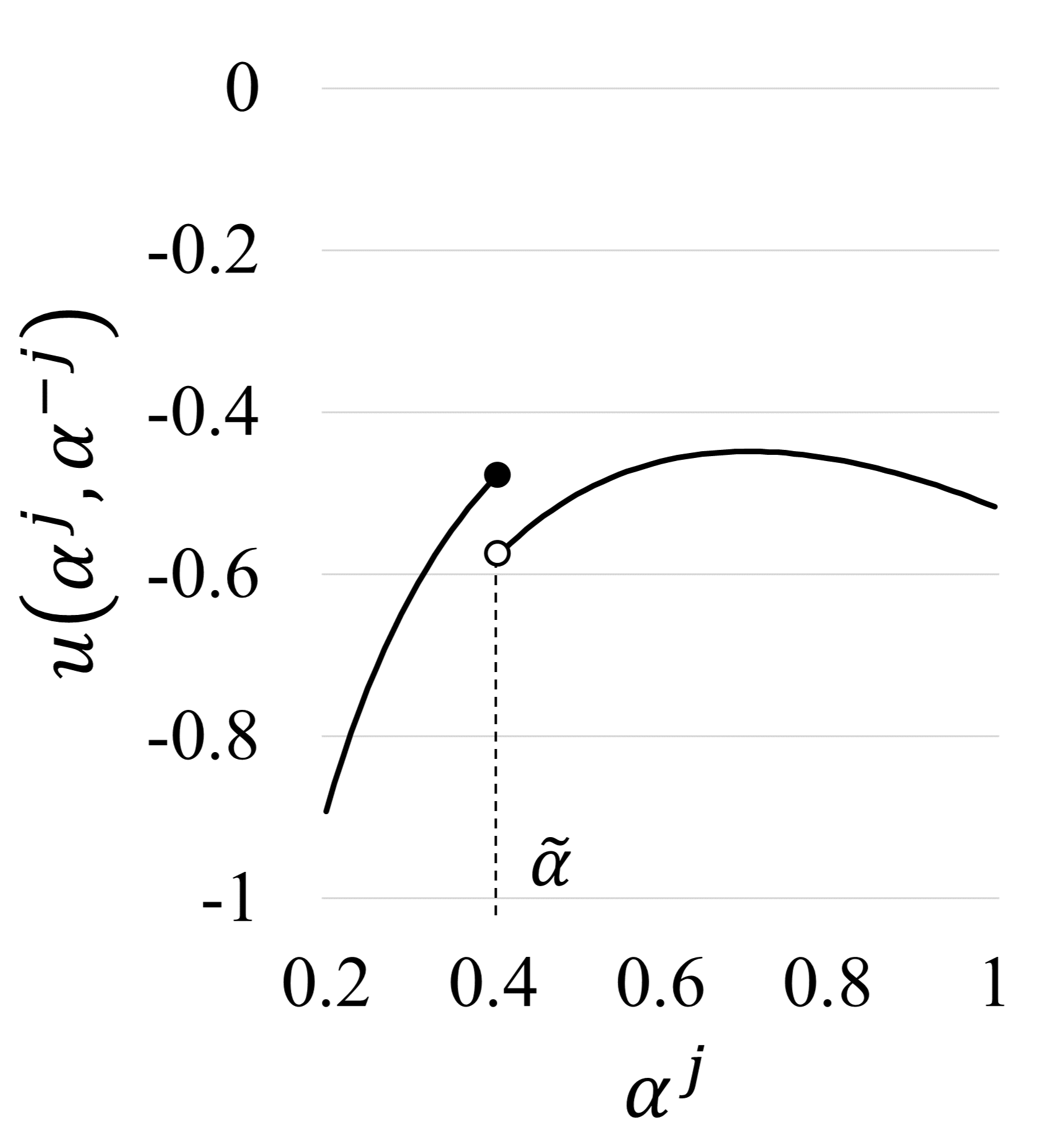}
            \caption[]%
            {}    
            \label{fig:BR1}
        \end{subfigure}
        \begin{subfigure}[b]{0.32\textwidth}  
            \centering 
            \includegraphics[width=\textwidth]{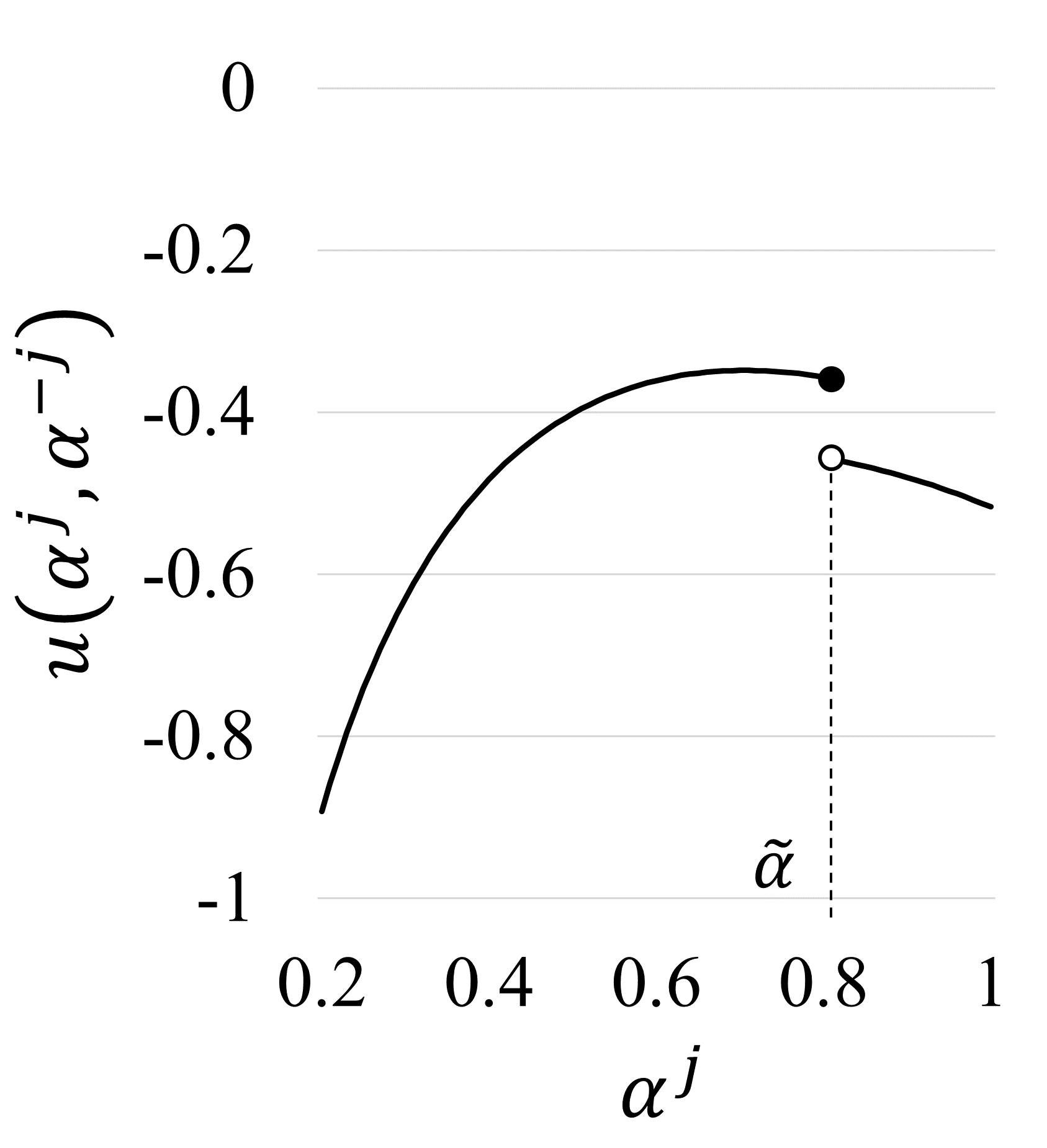}
            \caption[]%
            {}    
            \label{fig:BR2}
        \end{subfigure}
        \begin{subfigure}[b]{0.32\textwidth}   
            \centering 
            \includegraphics[width=\textwidth]{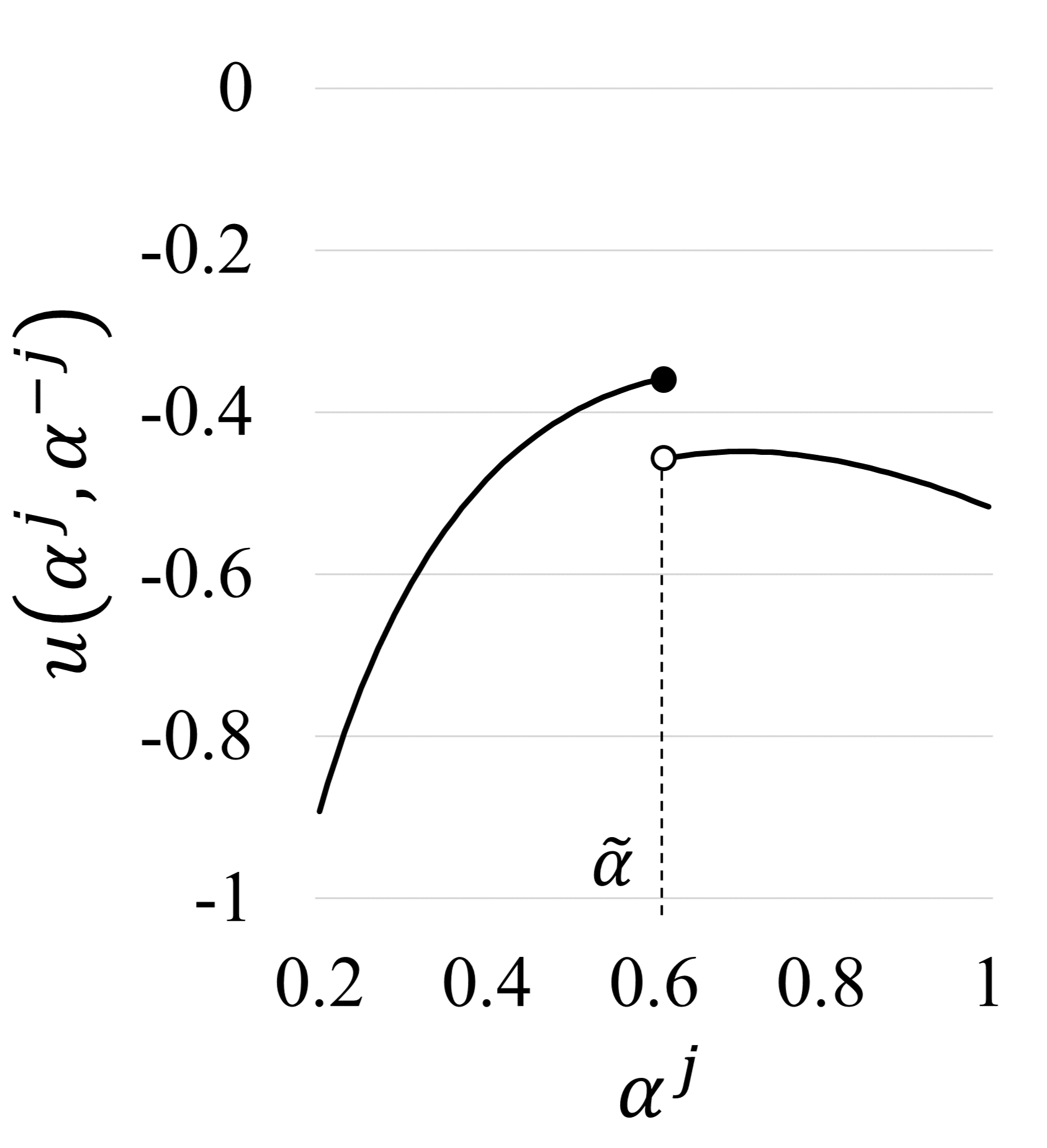}
            \caption[]%
            {}    
            \label{fig:BR3}
        \end{subfigure}
\centering
\caption[]{Illustration of the best response}
\label{fig:BR}
\end{figure}

To calculate Nash equilibrium for the game, we begin by rewriting the utility function of each player $j$ as follows:
\begin{equation}
\label{eq:utilityfcn}
    u(\alpha^j,\alpha^{-j})=\begin{cases}
        (\log\alpha^j-\alpha^j+1)-c-\alpha^j \kappa\tau\frac{\alpha^{-j} \hat i}{\alpha^{-j} \hat s + \alpha^{-j} \hat i + \alpha^{*R} \hat r}d, & \text{$\alpha^j>\tilde{\alpha}$}\\
        (\log\alpha^j-\alpha^j+1)-\alpha^j \kappa\tau\frac{\alpha^{-j} \hat i}{\alpha^{-j} \hat s + \alpha^{-j} \hat i + \alpha^{*R} \hat r}d, & \text{$\alpha^j \leq \tilde{\alpha}$}
        \end{cases}
\end{equation}
\noindent We take its first and second derivatives while holding $j$ constant, and we obtain
\textcolor{black}{
\begin{align}
    \frac{\partial_-}{\partial\alpha^j}u(\alpha^j,\alpha^{-j}) &=\frac{1}{\alpha^j}-1-\kappa\tau\frac{\alpha^{-j} \hat i}{\alpha^{-j} \hat s + \alpha^{-j} \hat i + \alpha^{*R} \hat r}d,  \label{eq:u'}\\
    \frac{\partial^2_-}{(\partial\alpha^j)^2}u(\alpha^j,\alpha^{-j}) &=-\frac{1}{(\alpha^j)^2}.   \label{eq:u''}
\end{align}
}
\textcolor{black}{Note we take the left derivative to take into account the discontinuity at $\tilde{\alpha}$ in the utility function.} The second derivative of $u(\alpha^j,\alpha^{-j})$, shown in Eq.\eqref{eq:u''}, is negative for all $\alpha^{j}$ in $ (0,1]$, indicating that $u(\alpha^j,\alpha^{-j})$ is concave for any fixed $\alpha^{-j}$.  
Thus, with the discontinuity at $\alpha^j=\tilde{\alpha}$, the best response of player $j$ with respect to $\alpha^{-j}$, which maximizes the utility of the agent $j$, is either the point where the first derivative is equal to 0 (Figure \ref{fig:BR}(a),(b)) or the discontinuity point $\tilde{\alpha}$ (Figure \ref{fig:BR}(c)). 

Next, let us define a function $g^*$:
\begin{definition}
\label{def:g*}
    $g^*(\alpha^{-j})$ is the value of $\alpha^j$ such that $\frac{\partial_-}{\partial \alpha^j}u(\alpha^j,\alpha^{-j})=0$.
\end{definition}


\noindent The function \( g^* \) represents the vertex of the utility function \( u \) with respect to \( \alpha^j \), depending on \( \alpha^{-j} \). In our context, it denotes the utility-maximizing action for individual \( j \) when all other agents adhere to \( \alpha^{-j} \). \textcolor{black}{To be precise, let us emphasize that whether \( g^*(\alpha^{-j}) \) is indeed the utility-maximizing action depends on the magnitude of discontinuity at $\tilde{\alpha}$. This aspect will be elaborated on further in subsequent discussions.}

We can rewrite $g^*(\alpha^{-j})$ in a closed form from Eq.\eqref{eq:u'}, as follows:
\begin{equation}
\label{eq:g*}
    g^*(\alpha^{-j})=\frac{\alpha^{-j} \hat s + \alpha^{-j} \hat i + \alpha^{*R} \hat r}{\alpha^{-j} \hat s + \alpha^{-j} \hat i + \alpha^{*R} \hat r + \kappa \tau d \alpha^{-j} \hat i}.
\end{equation}
In Eq.\eqref{eq:g*}, both the numerator and denominator are non-negative, with the denominator always equal to or greater than the numerator, ensuring that the value of \(g^*(\alpha^{-j})\) falls within the interval \((0,1]\). Additionally, \(g^*\) is a decreasing function with respect to \( \alpha^{-j} \). Intuitively, as the social activity of others increases, so does the infection risk for individual \(j\), prompting a reduction in \(j\)'s social activity level in response. 

Using the function $g^*$, we  define a variable $\Bar{\alpha}$:

\color{black}\begin{definition}
\label{def:a-bar}
$\bar{\alpha}$ is $\alpha$ such that $g^*(\alpha)=\alpha$. 
\end{definition}

\noindent The condition $g^*(\alpha)= \alpha$ describes a situation where, the utility function \(u\) peaks  -- its first derivative is zero -- at \(\alpha\) when all other agents follow \(\alpha\). With this definition, we can easily obtain $\bar{\alpha}$ in a closed form by solving Eq.\eqref{eq:g*}: 

\begin{equation}
\label{eq:a-bar}
    \bar{\alpha}=\frac{\hat s + \hat i - \alpha^{*R} \hat r + \sqrt{(\hat s + \hat i - \alpha^{*R} \hat r)^2+4\alpha^{*R} \hat r (\hat s + \hat i + \kappa \pi d \hat i)}}{2 (\hat s + \hat i + \kappa \pi d \hat i)}.
\end{equation}


As remarked in Definition 1, \(\bar{\alpha}\) may not represent the utility-maximizing action due to the discontinuity in the utility function. If it indeed is the utility-maximizing action, then \(\bar{\alpha}\) represents the optimal decision for individual agents; their best response is to follow the action chosen by others. This is the very definition of a Nash equilibrium, where \(\bar{\alpha}\) becomes the action universally adopted in equilibrium. 

While \(\bar{\alpha}\) may not be the utility-maximizing action depending on the magnitude of the discontinuity in $u$, it always is if there is no discontinuity. Thus, when there is no noncompliance penalty ($c = 0$), \(\bar{\alpha}\) maximizes the utility, hence  is the Nash equilibrium action.  In other words, \(\bar{\alpha}\) can be interpreted as the Nash equilibrium activity level in the absence of a social distancing policy.  This voluntary level of social activity is the result of balancing the benefits of social interaction against the risks of disease transmission. 

As discussed above, the presence of a discontinuity in the utility function $u$ makes this analysis complicated, and a closer examination is needed to identify the Nash equilibrium action. Before presenting Nash equilibrium solutions, let us establish the best response for individual $j$ in this game through the following proposition.

\begin{proposition}
\label{thm:BR}
The best response of individual $j$ for a given strategy of other players, $\alpha^{-j}$ is
\begin{equation*}
    BR(\alpha^{-j}) = \begin{cases}
                        \tilde{\alpha} \text{ with probability } p\\
                        g^*( \alpha^{-j}) \text{ with probability } 1-p\\
                      \end{cases}
\end{equation*}
for some $p \in [0,1]$. 
\end{proposition}
\noindent Proof of the proposition is in Appendix A. Proposition \ref{thm:BR}  states that the best response for individual $j$ is either to follow the policy ($\tilde{\alpha}$), or follow the action chosen by others. This proposition guarantees the optimality of the Nash equilibrium that we calculate hereafter. 

\color{black}

\begin{figure*}
        \centering
        \begin{subfigure}[b]{0.475\textwidth}
            \centering
            \includegraphics[width=\textwidth]{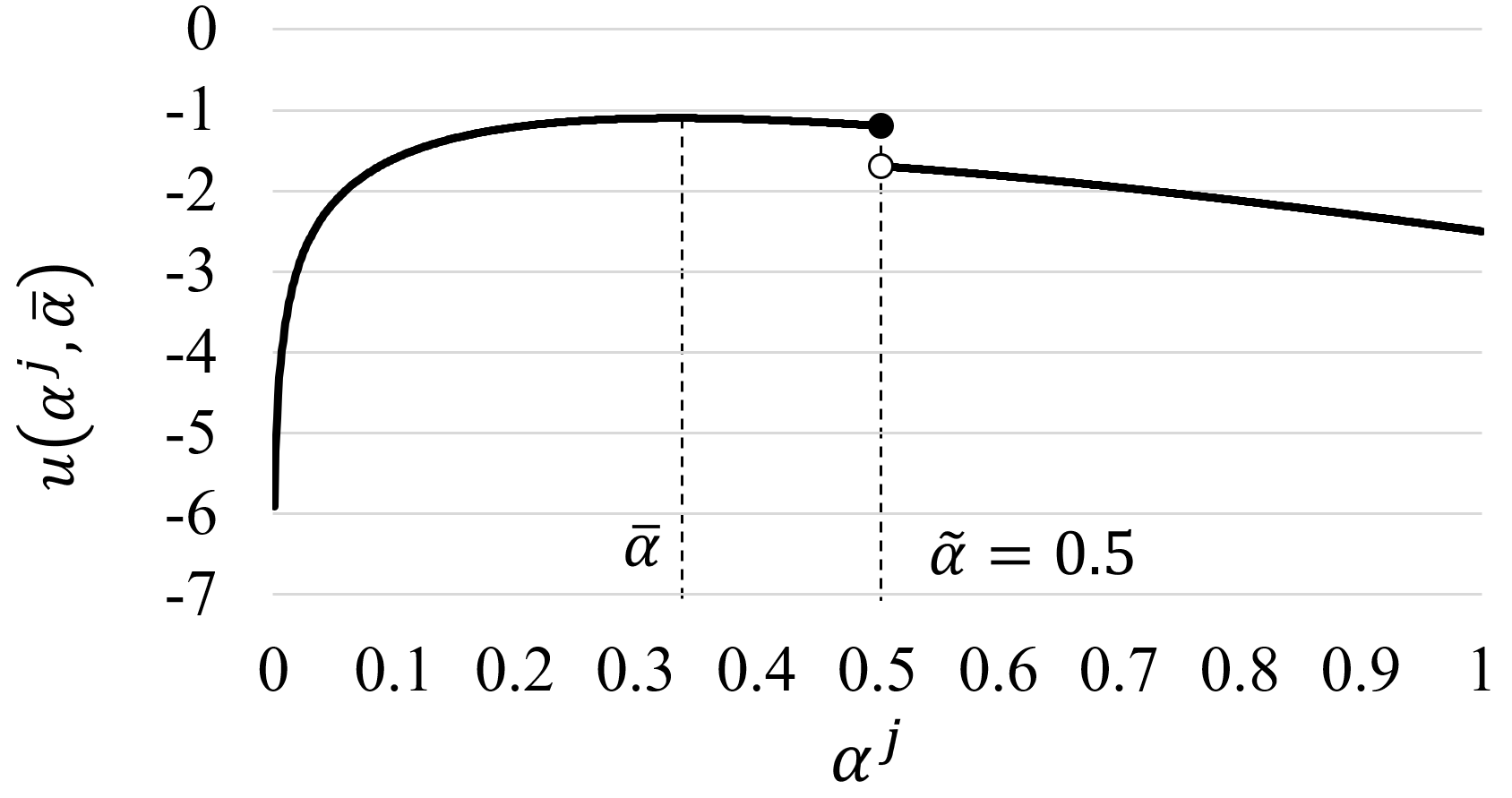}
            \caption[]%
            {}    
            \label{fig:NE1}
        \end{subfigure}
        \hfill
        \begin{subfigure}[b]{0.475\textwidth}  
            \centering 
            \includegraphics[width=\textwidth]{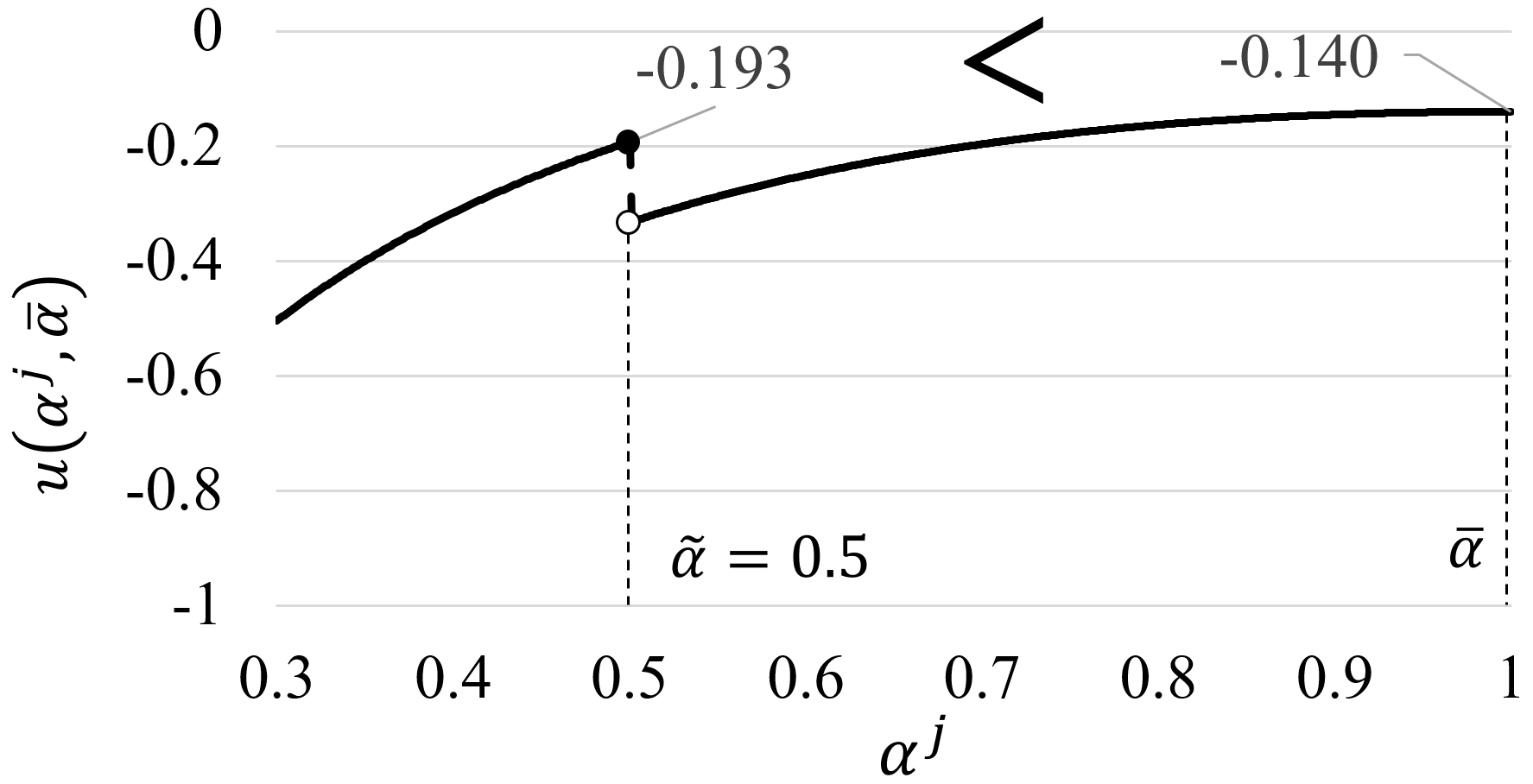}
            \caption[]%
            {}    
            \label{fig:NE2}
        \end{subfigure}
        \vskip\baselineskip
        \begin{subfigure}[b]{0.475\textwidth}   
            \centering 
            \includegraphics[width=\textwidth]{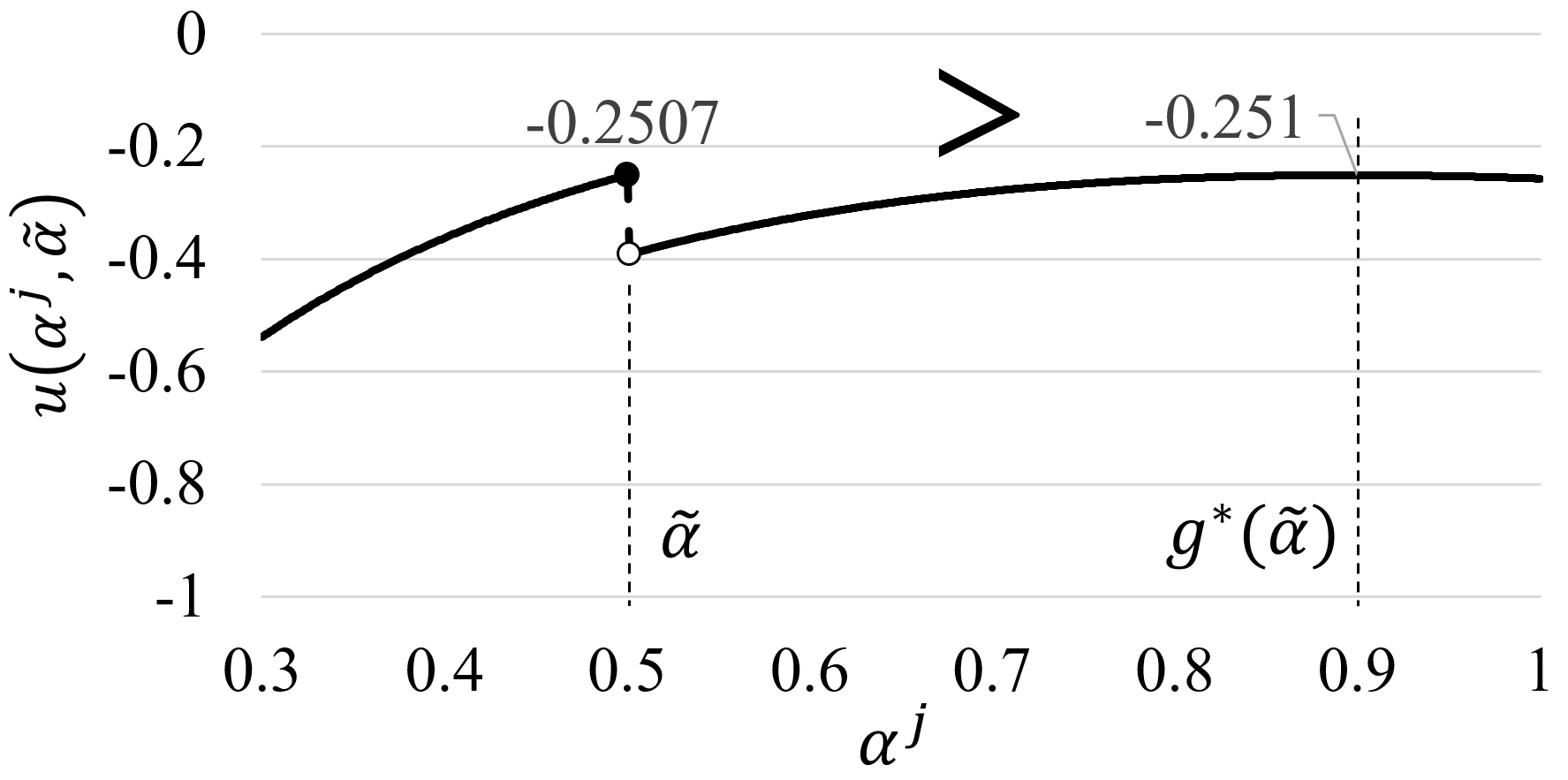}
            \caption[]%
            {}    
            \label{fig:NE3}
        \end{subfigure}
        \hfill
        \begin{subfigure}[b]{0.475\textwidth}   
            \centering 
            \includegraphics[width=\textwidth]{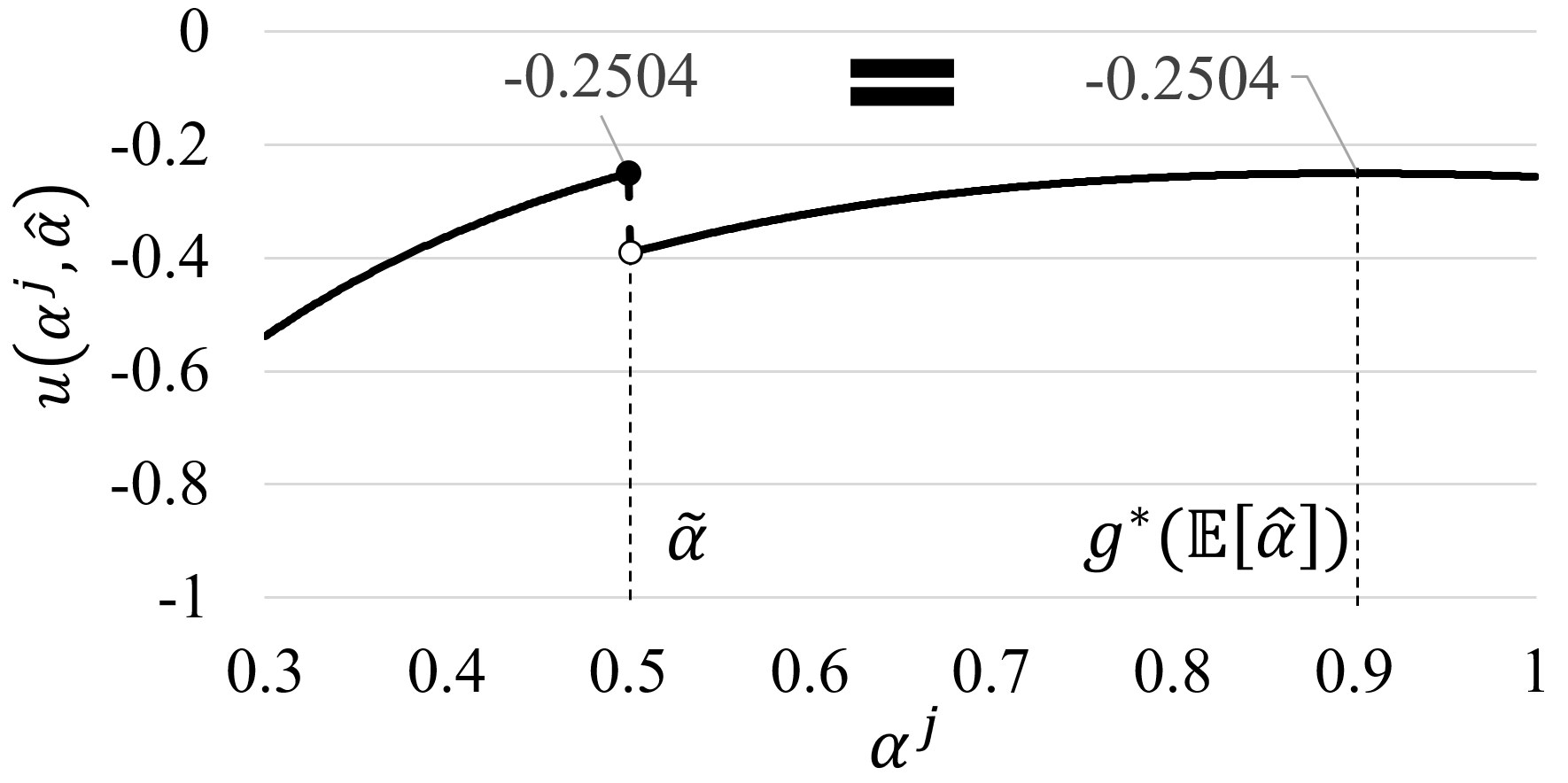}
            \caption[]%
            {}    
            \label{fig:NE4}
        \end{subfigure}
        \caption[ Examples of a Nash equilibrium ]
        {Examples of a Nash equilibrium} 
        \label{fig:NE-examples}
    \end{figure*}

Figure \ref{fig:NE-examples} presents the graphical illustration of the four Nash equilibrium cases, and we provide a detailed discussion below. Recall that at $\Bar{\alpha}$, the first derivative of $u$ is zero. When $\Bar{\alpha}$ is no greater than $\tilde{\alpha}$ ($\bar{\alpha} \leq \tilde{\alpha}$), the apex of $u$ is located to the left of the discontinuity drop, and this condition ($\frac{\partial_- u}{\partial \alpha^{h}} = 0$) suffices for $\Bar{\alpha}$ to represent a pure Nash equilibrium. This case is shown in Figure \ref{fig:NE-examples}(a). In contrast, when $\Bar{\alpha}$ is greater than $\tilde{\alpha}$ ($\bar{\alpha} > \tilde{\alpha}$), we have two cases to \ consider. The first case, illustrated in Figure \ref{fig:NE-examples}(b), is when the utility function for individual $j$ is greater at $\alpha^j = \Bar{\alpha}$ than at $\Tilde{\alpha}$:  $u(\Bar{\alpha},\Bar{\alpha}) \geq u(\tilde{\alpha},\Bar{\alpha})$. In this case, $\Bar{\alpha}$ is still the best response for agent $j$, and thus $\Bar{\alpha}$ is again a pure Nash equilibrium action. The second case is when $u(\Bar{\alpha},\Bar{\alpha}) <  u(\tilde{\alpha},\Bar{\alpha})$. This happens when the drop at $\Tilde{\alpha}$ is so significant that the apex of $u$ is lower than the value of $u$ at $\Tilde{\alpha}$. Since it yields a higher utility value, $\tilde{\alpha}$ will be a better response for individual $j$ than $\Bar{\alpha}$. Thus $\Bar{\alpha}$ is no longer a Nash equilibrium point, and we need to find a new candidate for Nash equilibrium.

To do so, we reconstruct the utility function for the case where $\alpha^{-j} = \Tilde{\alpha}$. It turns out that for this new $u(\alpha^j,\Tilde{\alpha})$, its apex $g^*(\tilde{\alpha})$ always exists to the right of 

To do so, we reconstruct the utility function for the case where $\alpha^{-j} = \Tilde{\alpha}$. It turns out that for this new $u(\alpha^j,\Tilde{\alpha})$, its apex $g^*(\tilde{\alpha})$ always exists to the right of the drop $\tilde{\alpha}$. This is because $g^*$ is a decreasing function of $\alpha^{-j}$ and we have $\tilde{\alpha} < \bar{\alpha}$, leading to $g^*(\tilde{\alpha}) > g^*(\Bar{\alpha}) = \Bar{\alpha} > \tilde{\alpha}$. Again, we have two cases to consider. In the first case when $u(\tilde{\alpha},\tilde{\alpha}) \geq  u(g^*(\tilde{\alpha}),\tilde{\alpha})$ as in Figure \ref{fig:NE-examples}(c),  the best response by individual $j$ is $\tilde{\alpha}$ and thus $\tilde{\alpha}$ is a Nash equilibrium point. The other case is when $u(\tilde{\alpha},\tilde{\alpha}) < u(g^*(\tilde{\alpha}),\tilde{\alpha})$, and  $\tilde{\alpha}$ is not the best response for agent $j$. In this case, a pure Nash equilibrium does not exist and we must instead seek a mixed Nash equilibrium. \textcolor{black}{We can prove the existence of a mixed Nash equilibrium in our game, and a proposition and its proof is provided in Appendix A. Below, we briefly present the result.}

\color{black}A mixed Nash equilibrium solution is defined by a strategy $\hat{\alpha}$ in the following equations:
 
\begin{align}
    &\hat{\alpha}=  \begin{cases}
                    \tilde{\alpha} & \text{with probability } p\\
                    g^*( \mathbb{E} \left[ \hat{\alpha} \right]) & \text{with probability } 1-p 
                    \end{cases} \label{eq:mu-hat-def} \\ 
    &u(\tilde{\alpha},\mathbb{E} \left[ \hat{\alpha} \right]) =u(g^*(\mathbb{E} \left[ \hat{\alpha} \right]),\mathbb{E} \left[ \hat{\alpha} \right] ). \label{eq:mu-hat-cond}
\end{align}

\noindent Eq.\eqref{eq:mu-hat-def} represents an action $\hat{\alpha}$ by an individual under which it chooses $\Tilde{\alpha}$ with probability $p$ and $g^*(\mathbb{E} \left[ \hat{\alpha} \right])$ with probability $(1-p)$. With Eq.\eqref{eq:mu-hat-def} and \eqref{eq:mu-hat-cond}, we seek to obtain the value of $\mathbb{E} \left[ \hat{\alpha} \right]$ and $p \in (0,1)$ such that both $\tilde{\alpha}$ and $g^*(\mathbb{E} \left[ \hat{\alpha} \right])$ are the best responses by player $j$ when all other individuals behave according to $\hat{\alpha}$. If such $\hat{\alpha}$ exists, $(p,1-p)$ for action $(\tilde{\alpha},g^*(\mathbb{E} \left[ \hat{\alpha} \right]))$ is the best response for individual $j$ and is a Nash equilibrium solution. An example of such solution is shown in Figure \ref{fig:NE-examples}(d). 




A complete set of Nash equilibrium solution for our game is summarized in Table \ref{tab:NE}. \color{black} We use these Nash equilibrium solutions in simulation to obtain results for numerical experiments. For the completeness of analysis, we need to ensure that the cases shown in Table \ref{tab:NE} do not overlap and that the Nash equilibrium in each case is unique. Both of those hold true in our game, and we provide its proof in the Appendix A. 

\begin{table}[htb]
\caption{Nash equilibrium for each case}
\label{tab:NE}
\begin{tabular} {clc}
\hline
Case & Condition & Nash Equilibrium \\
\hline
1&$\bar{\alpha}\leq\tilde{\alpha}$&$\bar{\alpha}$\\
2&$\bar{\alpha}>\tilde{\alpha}$ and $u(\bar{\alpha},\bar{\alpha}) \geq u(\tilde{\alpha},\bar{\alpha})$ & $\bar{\alpha}$ \\
3&$\bar{\alpha}>\tilde{\alpha}$ and $u(\tilde{\alpha},\tilde{\alpha}) \geq u(g^*(\tilde{\alpha}),\tilde{\alpha})$ & $\tilde{\alpha}$ \\
\multirow{2}{*}{4} & $\bar{\alpha}>\tilde{\alpha}$ and & \multirow{2}{*}{$\hat{\alpha}$}\\
 & $\quad u(\bar{\alpha},\bar{\alpha}) < u(\tilde{\alpha},\bar{\alpha})$ and $u(\tilde{\alpha},\tilde{\alpha}) < u(g^*(\tilde{\alpha}),\tilde{\alpha})$ & \\
\hline
\end{tabular}
\end{table}

\section{Experimental Results}
\label{sec:result}

\subsection{Experimental settings and overview of Nash equilibrium solutions}
\label{subsec:result-setting}

To conduct our numerical experiments, we set the model parameter values based on the information from the COVID-19 case, and for some parameters for which no relevant data is available, we make assumptions. These values are summarized in Table \ref{tab:parameter}. For the basic reproductive number, we use $R_0 = 3.1$ which is an estimate from the data in Korea \citep{jeong2020estimation-R0}. This estimate is within the range of  other estimations reported in the literature, e.g. \cite{alimohamadi2020estimate-R0}.  \textcolor{black}{For the epidemic model parameters, $\gamma_1$ and $\gamma_2$ in Eq.\eqref{eq:SIR}, we refer to the literature for related information such as mean incubation period \citep{alene2021serial-incubation,mcaloon2020incubation-incubation, quesada2021incubation-incubation} and quarantine period \citep{chtourou2020staying-quarantine,CDC-quarantine}, and make reasonable assumptions on these parameters. For $\beta$, we compute its value by using $R_0 = \beta / \gamma_1$.}

\begin{table}[hbt]
    \caption{Model parameters}
    \centering
    \begin{tabular}{lc}
        \hline
        Parameter & Value \\
        \hline
        basic reproductive number, $R_0$ &  3.1 \\
        \textcolor{black}{(mean incubation period)$^{-1}$,  $\gamma_1$ } & 1/7 \\
        (quarantine period)$^{-1}, \gamma_2$  & 1/14 \\
        \textcolor{black}{transmission rate parameter}, $\beta$ & 3.1/7 \\
        \textcolor{black}{sensitivity to  infection probability}, $d$ & 1858 \\
        \hline
    \end{tabular}
    \label{tab:parameter}
\end{table}

For parameter $d$, used in Eq.\eqref{eq:obj-SI} to capture the perceived cost of getting infected, we estimate its value by utilizing the data from the COVID-19 case in Korea. Specifically, we collected data for the period of February 24th to March 1st of 2020, which is the very early stage of the disease spread. As there was no social distancing policy (i.e., $\tilde{\alpha}_t=1$ and $c_t = 0$) during this period, we assume the level of social activity maintained by the public is the collective result of individual decision making according to Eq.\eqref{eq:a-bar}. In Eq.\eqref{eq:a-bar}, $\bar{\alpha}$ denotes the level of social activity by the public (susceptible and infected) and $\alpha^{*R}$ for the recovered group. We estimate $\bar{\alpha}$ by comparing average daily mobility data\footnote{Daily mobility data based on mobile phone traffic is provided by a major telecommunication company in Korea, which can be accessed at https://data.kostat.go.kr/social/moblilePopMoveInfoPage.do.} before COVID-19 to that of February 24--March 1st period. For the recovered, $\alpha^{*R}$ is assumed to be 1 since it is the best response for the recovered group. Then, we use the daily count of confirmed COVID-19 cases to determine the epidemic state variables $(\hat{s},\hat{i},\hat{q}, \hat{r})$. These values are plugged into Eq.\eqref{eq:a-bar} to obtain an estimate for $d$.

With the model parameters set, we conduct experiments for various values of $c_t$ and $\tilde{\alpha}_t$ to explore the effect of a social distancing policy on individuals' compliance decisions and overall disease spread dynamics. A time unit for the numerical experiments is a day, and we run the simulation for the length of one year, $T = 365$. Note that for a simpler and interpretable analysis, we focus on scenarios where $c_t$ and $\tilde{\alpha}_t$ are fixed over $t$ -- i.e., a social distancing policy is set at the beginning and maintained throughout the time horizon.  

\begin{figure}[htb]
 \centering
        \begin{subfigure}[b]{0.32\textwidth}
            \centering
            \includegraphics[width=\textwidth]{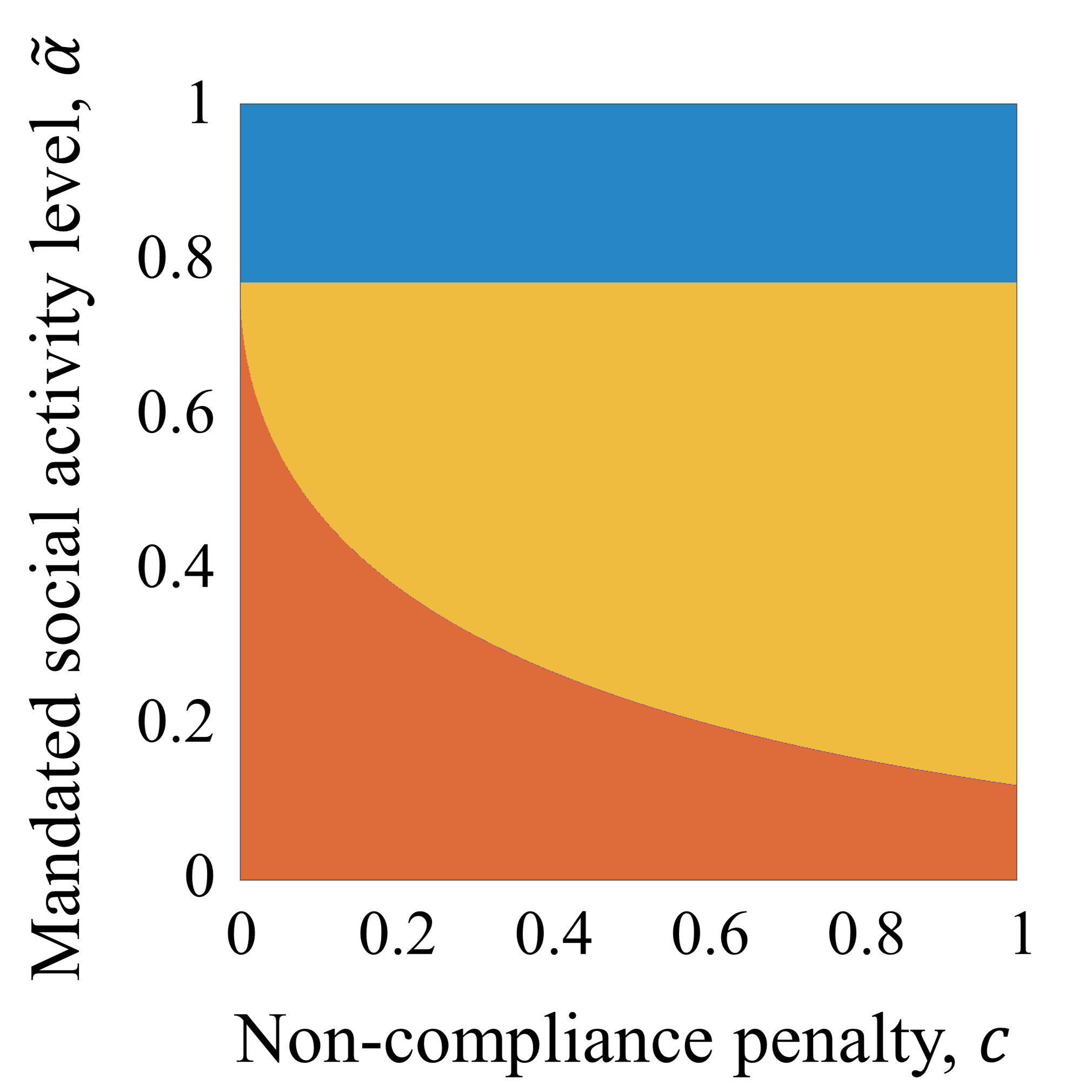}
            \caption[]%
            {}    
            \label{fig:cases1}
        \end{subfigure}
        \begin{subfigure}[b]{0.32\textwidth}  
            \centering 
            \includegraphics[width=\textwidth]{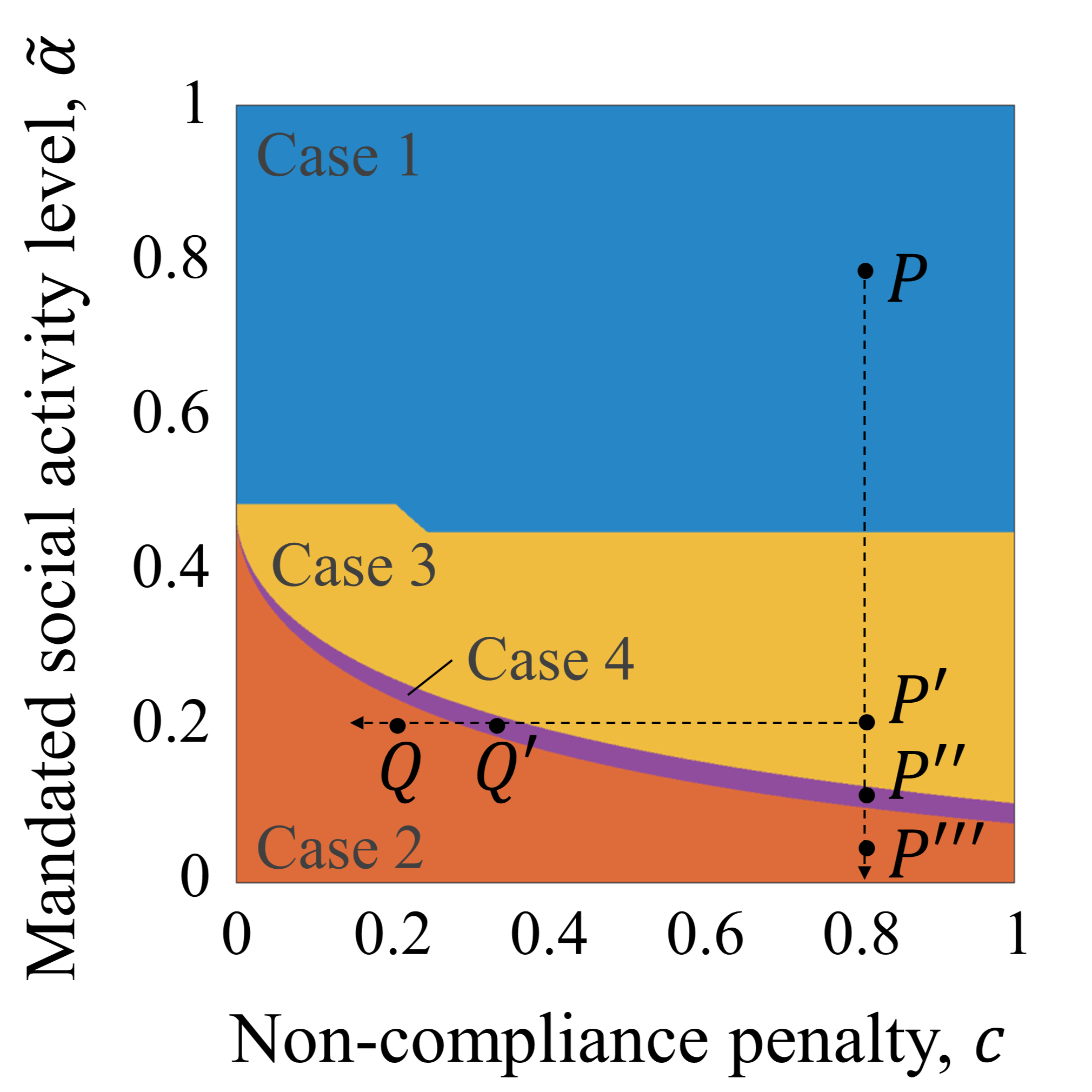}
            \caption[]%
            {}    
            \label{fig:cases2}
        \end{subfigure}
        \begin{subfigure}[b]{0.32\textwidth}   
            \centering 
            \includegraphics[width=\textwidth]{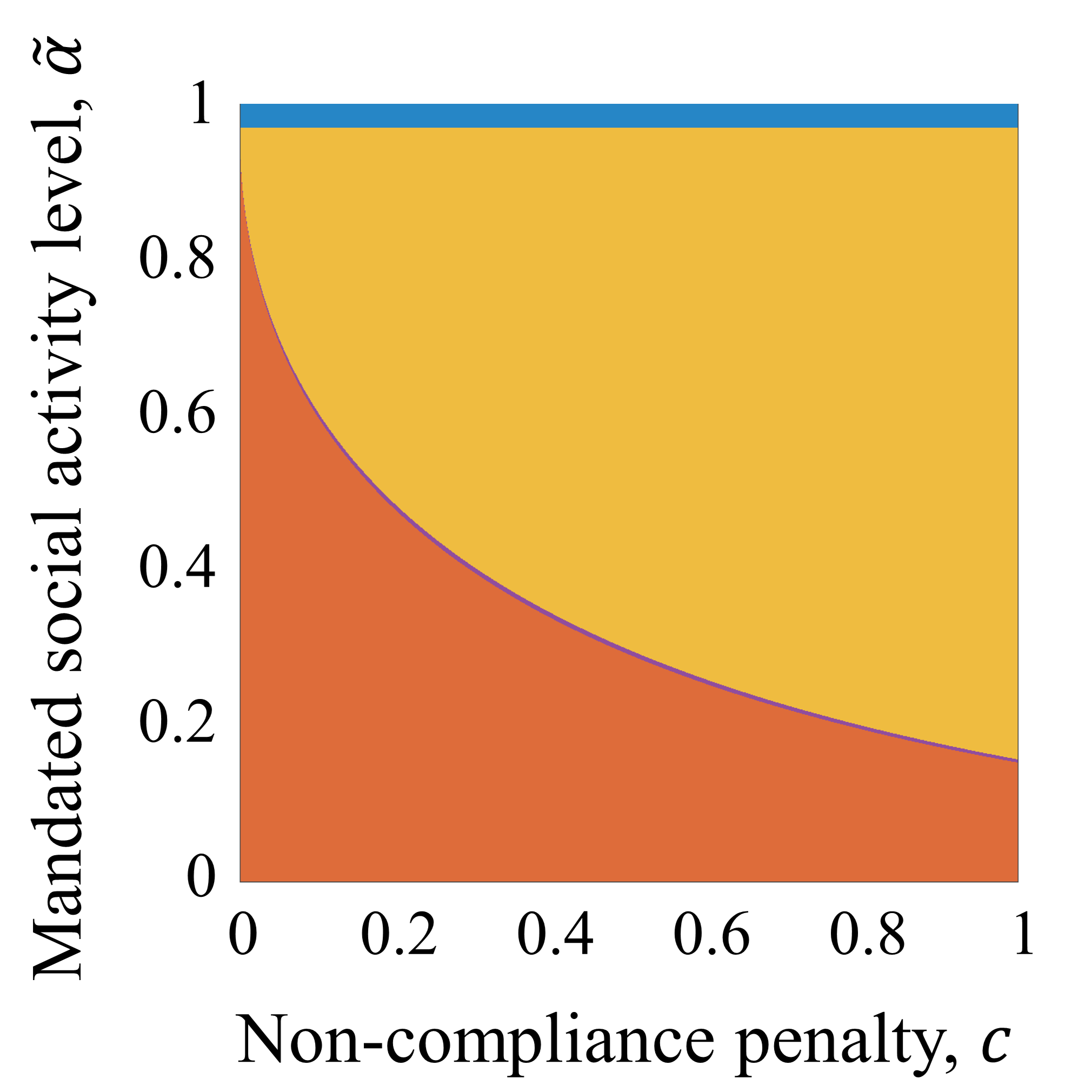}
            \caption[]%
            {}    
            \label{fig:cases3}
        \end{subfigure}
\centering
\caption[]{\textcolor{black}{Illustration of Nash equilibrium at three states in (a) early period:$(s,i,q,r)$ = (0.9431, 0.0356, 0.0168,	0.0045), (b) peak:(0.5885, 0.1233, 0.1656, 0.1226), and (c) post-peak period:(0.2235, 0.0033, 0.0171, 0.7561)}
}
\label{fig:NE_alpha}
\end{figure}

\color{black}
Figure \ref{fig:NE_alpha} shows plots of the Nash equilibrium solution of $\alpha$ on the ($\tilde{\alpha}_t, c_t$) space, which is a graphical presentation of Table \ref{tab:NE}. Since the equilibrium solution is determined for a specific system state, we  chose three states for plotting $\alpha_t$
.

\noindent In the figure, we use $c_t$ for x-axis and $\tilde{\alpha}_t$ for y-axis for better presentation. In Figure \ref{fig:NE_alpha}(b), we show how the equilibrium solution $\alpha_t^*$ changes within the $(c_t,\tilde{\alpha}_t)$ space. 
Let us start with point $P = (0.8, 0.8)$ in the region ``Case 1'', where the Nash equilibrium solution for $\alpha_t$ is $\bar{\alpha}_t$ when the mandate from the social distancing policy is $\tilde{\alpha}_t$. Note from Table \ref{tab:NE} that $\tilde{\alpha}_t \geq \bar{\alpha}_t$ in this region. Thus, this Nash equilibrium means that when the required social distancing is too loose, even weaker than the voluntary social distancing level practiced  for self-protection, individuals behave conservatively at the level of voluntary social distancing. Moving down along the vertical line ($\overline{PP'''}$), we have $P' = (0.8, 0.2)$ in ``Case 3'' region. In this region, $\tilde{\alpha}_t < \bar{\alpha}_t$ and the Nash equilibrium solution becomes $\tilde{\alpha}_t$. The policy's social distancing mandate  is stricter than the voluntary social distancing, and  individuals are willing to follow the mandated level of social distancing. Further down at $P'' = (0.8, 0.12)$, residing in region ``Case 4'' where the social distancing requirement becomes even stricter, then the Nash equilibrium for individuals social distancing level is $\hat{\alpha}_t$, a mixed Nash equilibrium; they follow the social distancing policy only probabilistically. Then we have $P''' = (0.8, 0.05)$, where the required level $\tilde{\alpha}_t$ is very low. In this region ``Case 2'', the Nash equilibrium solution is $\bar{\alpha}'_t$, and $\bar{\alpha}'_t > \tilde{\alpha}_t$. That is, the social distancing requirement is too strict to follow, and people resort back to their natural, voluntary social distancing. 

Now let us examine along the horizontal line ($\overline{P'Q}$). At $P'$ in ``Case 3'', the Nash equilibrium solution is $\tilde{\alpha}_t$. Moving leftward along the line as the non-compliance penalty is reduced, we have $Q'$ in ``Case 4'' where the solution is now the mixed Nash equilibrium $\hat{\alpha}_t$. That is, at the same level of social distancing requirement, individuals start to deviate from full compliance when the noncompliance penalty drops below certain level. If the penalty is reduced further, we have $Q$ in ``Case 2'' and the Nash equilibrium solution is $\bar{\alpha}'_t$; the penalty is too low to enforce compliance, and individuals' social distancing practice is driven by the voluntary, self-protection behavior.  

\subsection{Intensity of social distancing policy}
\label{subsec:result-intensity}

\color{black}
We investigate  impacts of the intensity of social distancing on the spread of the epidemic disease. In Figure \ref{fig:w,wo-policy}, we present a few example outputs\footnote{\textcolor{black}{To generate various patterns of $i(t)$ curves, we use different value for $d$ than shown in Table \ref{tab:parameter}. For all subsequent experiments, $d$ is set as specified in Table \ref{tab:parameter}.}}  from the experiments where we vary intensities of social distancing policy while keeping the penalty $c_t$ fixed. These results highlight the complex nature of how the transmission of the disease is influenced by $\tilde{\alpha}_t$.

\begin{figure}[htb]
\centering
\includegraphics[width=0.75\textwidth]{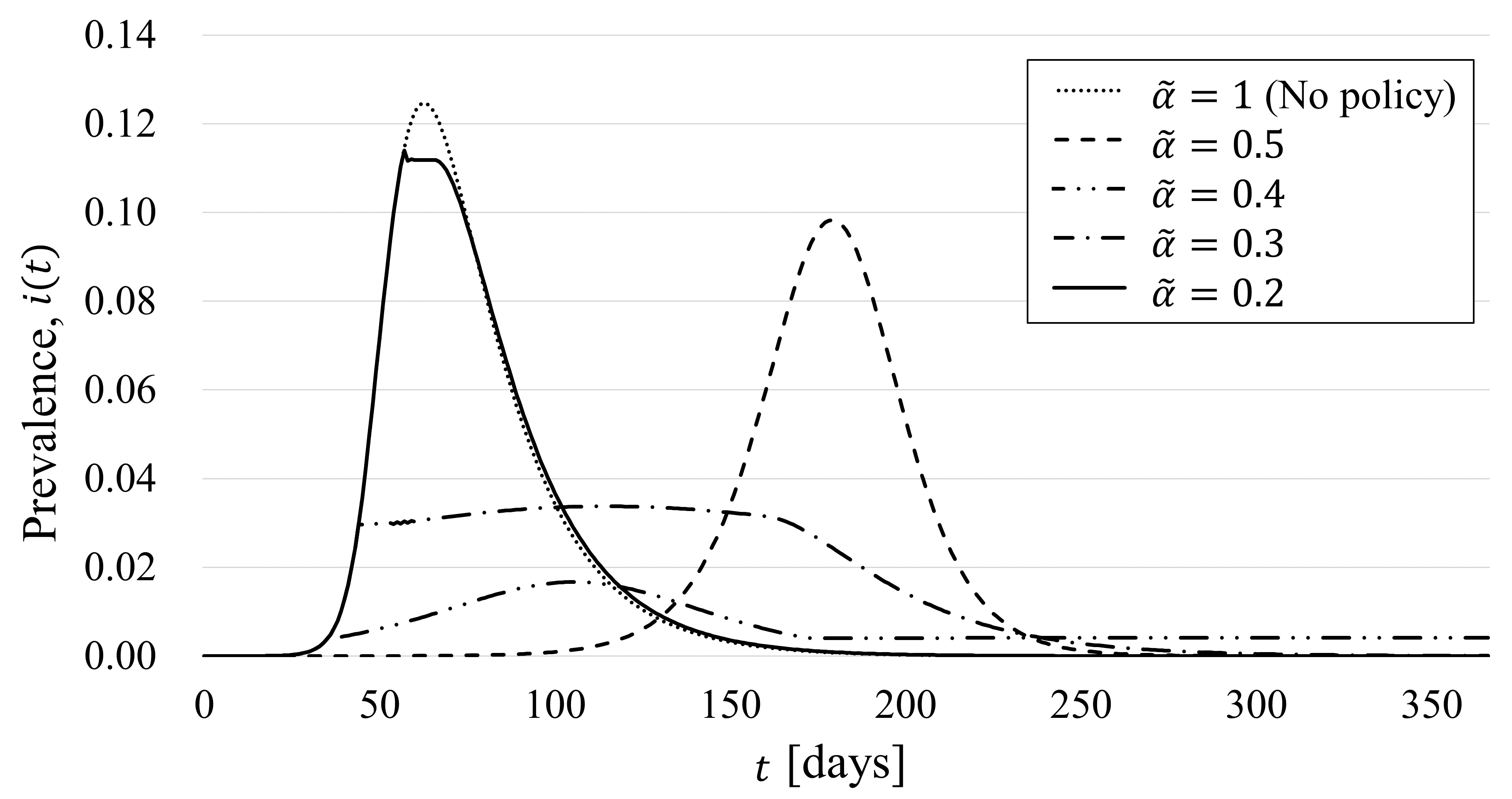}
\caption[]{\textcolor{black}{Prevalence under social distancing policy of various intensities, $\tilde{\alpha}_t$ }
}
\label{fig:w,wo-policy}
\end{figure}

Each curve shown in the figure is a solution to the proposed compartmental model -- Eq.\eqref{eq:SIR} and \eqref{eq:adjusted_dsdt} -- when individuals set their social activity level according to Eq.\eqref{eq:decision-R} and Table \ref{tab:NE}. It is observed that compared with no social distancing scenario, introducing a social distancing policy effectively mitigates the spread of the disease; reducing and delaying the peak of infections, and ultimately lowering the overall infection count. While such improvement is evidently expected, the most interesting observation from these curves is that the effectiveness from social distancing policies is not simply monotonically increasing as we raise its intensity; it is most effective when $\tilde{\alpha}_t$ is moderate, and becomes less effective when $\tilde{\alpha}_t$ is decreased below that level, eventually becoming almost ineffective when it is set too low. 

A particularly interesting result is when $\tilde{\alpha}_t = 0.2$, shown in a solid line. Under this very strict social distancing policy, prevalence $i(t)$ is almost identical to the no policy case ($\tilde{\alpha}_t = 1$). This seemingly counter intuitive result deserves some discussions, for which we revisit Figure \ref{fig:NE_alpha}. We observe that when $\tilde{\alpha}_t = 1$, the Nash equilibrium solution is  $\bar{\alpha}_t$ in all states, and for very low $\tilde{\alpha}_t$, the Nash equilibrium solution is $\bar{\alpha}'_t$. Thus, the difference between the no social distancing policy and very strict policy scenarios depends on the discrepancy between $\bar{\alpha}_t$ and $\bar{\alpha}'_t$. According to Eq.\eqref{eq:a-bar}, with other variables being equal, the difference between $\bar{\alpha}_t$ and $\bar{\alpha}'_t$ is determined by $\alpha^{*R}_t$. The optimal action for the recovered population, $\alpha^{*R}_t$, is given by Eq.(7), and if the parameter $c_t$ is set such that $\alpha^{*R}_t$ is 1, then we have $\bar{\alpha}_t=\bar{\alpha}'_t$. This explains how  individuals’ level of social activities are the same between no social distancing policy ($\tilde{\alpha}_t=1)$  and near lockdown policy case ($\tilde{\alpha}_t=0.2)$, hence resulting in almost identical infection level. For generating results in Figure \ref{fig:w,wo-policy}, the non-compliance penalty $c_t$ was set to 0.35. If we used lager $c_t$ to make $\alpha^{*R}_t = \tilde{\alpha}_t$, $\bar{\alpha}'_t$ would be smaller; for much larger $c$, the equilibrium solution would then be found in ``Case 3'' to be $\tilde{\alpha}_t$ (Figure \ref{fig:NE_alpha}). Thus, this result highlights the importance of choosing proper level of $c_t$ when imposing a strict social distancing policy.

\begin{figure}[htb]
\centering
\includegraphics[width=0.9\textwidth]{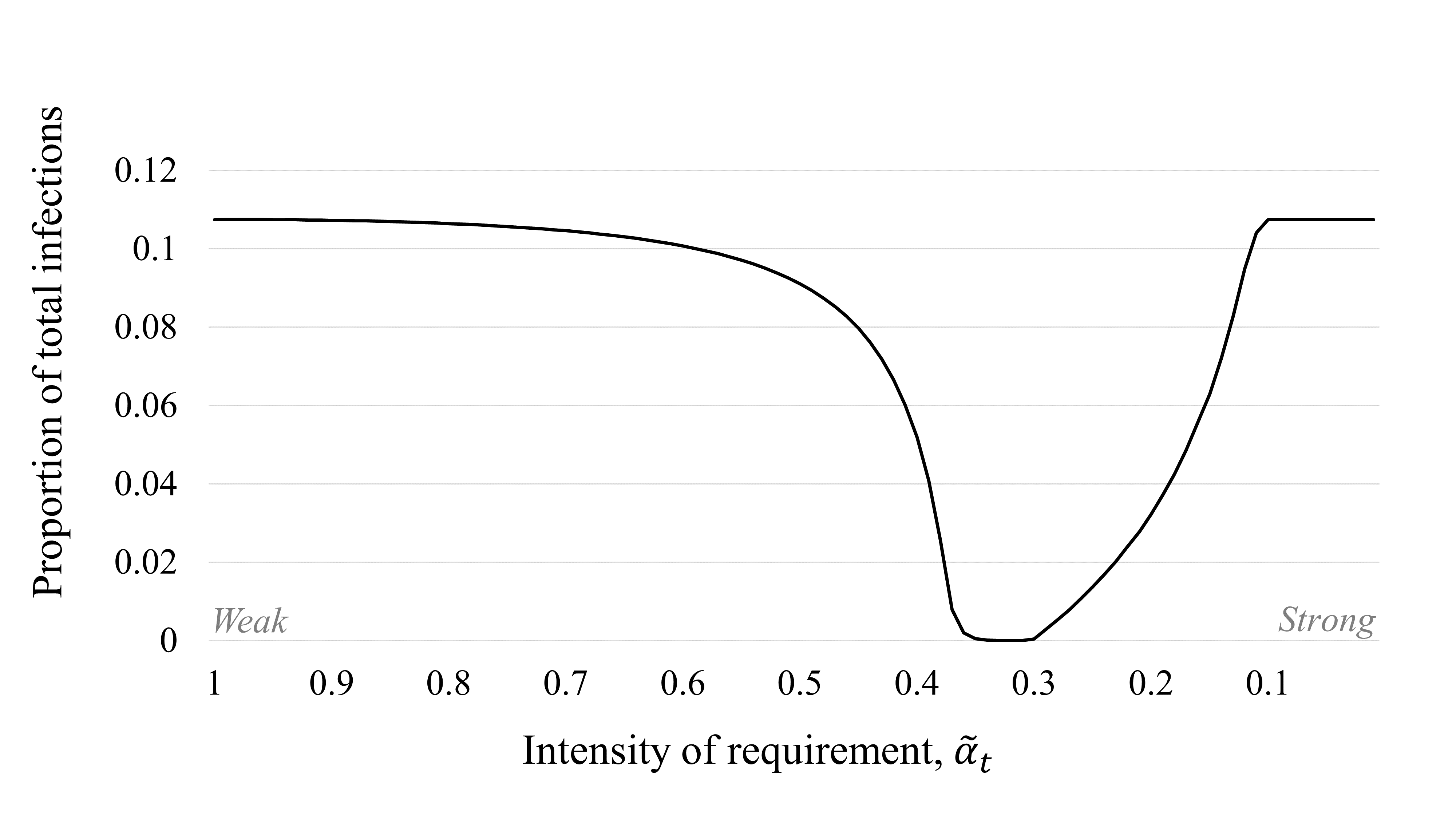} \caption{Proportion of the total infected population according to social distancing intensity with a fixed noncompliance penalty ($c_t=0.5$)}
\label{fig:alpha-tilde-result}
\end{figure}

In Figure \ref{fig:alpha-tilde-result}, we show  the total infection as a function of social distancing intensity, $\tilde{\alpha}_t$. Consistent with the results shown in Figure \ref{fig:w,wo-policy}, we observe a drop and rebound in the total number of infections as the intensity of the imposed policy becomes stronger. Initially, increasing the intensity of social distancing (i.e., decreasing $\tilde{\alpha}_t$) leads to a gradual decrease in total infections, and the decline becomes more pronounced when $\tilde{\alpha}_t$ falls below 0.5. However, when $\tilde{\alpha}_t$ reaches 0.3, the total infection starts to increase. This result is obtained consistently under various settings (Appendix B).

\begin{figure}[bht]
\centering
\includegraphics[width=0.75
\textwidth]{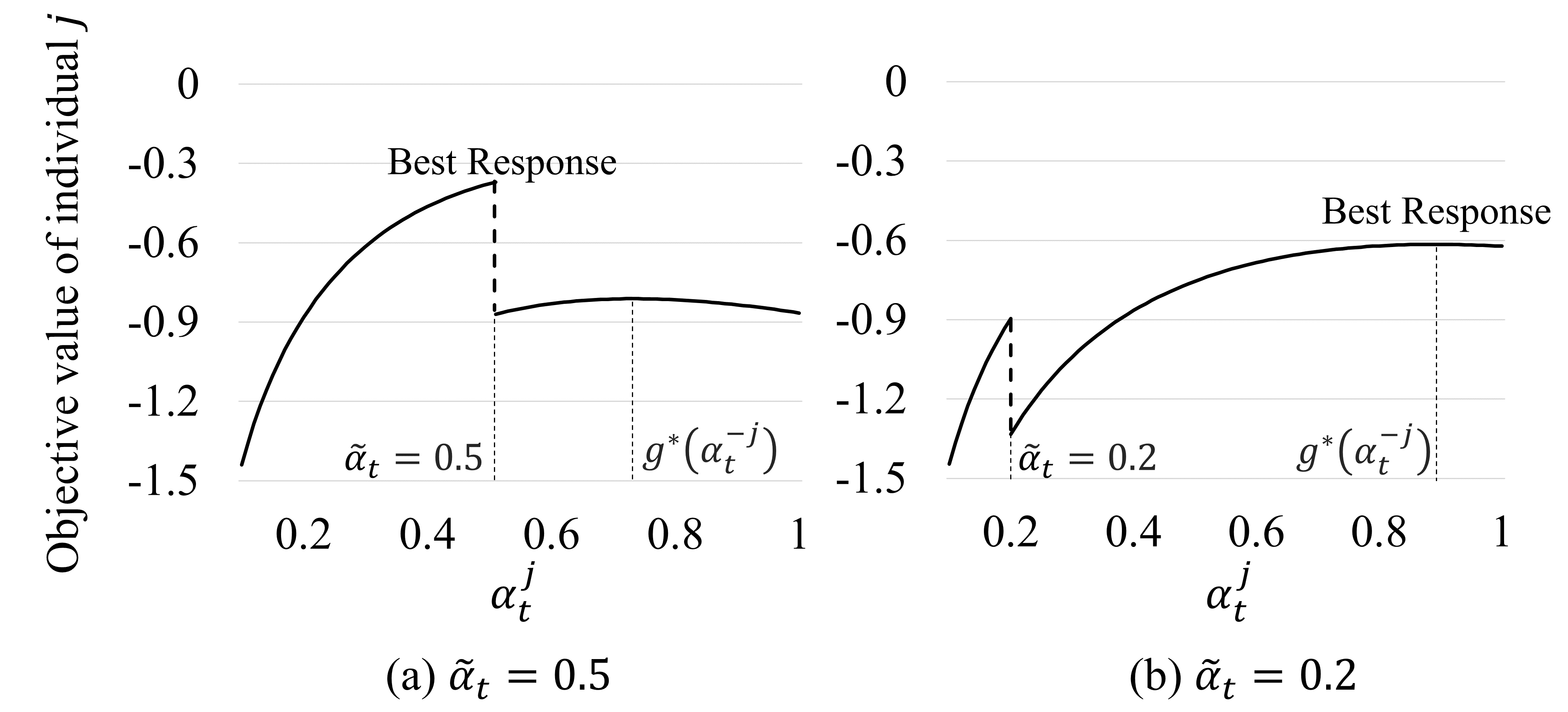}
\caption{Objective function as a function of social activity level of an individual $j$ given a fixed social activity level of other individuals: (a) $\tilde{\alpha}_t = 0.5$ (b) $\tilde{\alpha}_t = 0.2$} 
\label{fig:alpha-tilde-discussion}
{}
\end{figure}

Again, this rebound in the total infections can be attributed to the reduced degree of policy compliance from the community when the intensity of social distancing is too high relative to the level of epidemic spread. This can be better understood by referring to Figure \ref{fig:alpha-tilde-discussion}, which shows how the objective function changes as the social activity level of an individual $j$ varies. Note that these curves are drawn while keeping $\alpha_t^{-j}$ constant for illustration purposes. 

For a relatively weak social distancing intensity (e.g., $\tilde{\alpha}_t = 0.5$), individuals are motivated to comply with the policy as their best response is indeed $\tilde{\alpha}_t$, as shown in Figure \ref{fig:alpha-tilde-discussion}(a). On the contrary, when the social distancing intensity is strong, as illustrated in Figure \ref{fig:alpha-tilde-discussion}(b) for $\tilde{\alpha} = 0.2$, a higher objective value can be achieved by increasing their social activity level beyond $\tilde{\alpha}_t$, and this noncompliance option becomes their best response.

\color{black}

Results from these numerical experiments demonstrate that individuals may choose to disregard the social distancing policy when the intensity of social distancing is too strong relative to the degree of epidemic spread. These findings align with the empirical evidence in the literature, which suggests that the mathematical model and Nash equilibrium proposed in this study can serve as a theoretical framework to explain the observed behaviors regarding social distancing policies.

\subsection{Noncompliance penalty}
\label{subsec:result-penalty}

In the previous section, we demonstrated that an excessively high intensity of social distancing can decrease its effectiveness due to the non-compliance in the community. This suggests that the public health authority cannot be too aggressive in imposing social distancing policy. However, circumstances may require strong social distancing, for example, when the increase in total infections are getting out of control. In such cases, the authority may need to invoke a larger penalty in order to drive the community to comply  with a high-intensity social distancing policy. In this section, we investigate the effect of non-compliance penalty on the spread of disease.  

\begin{figure}[hbt!]
\centering
\includegraphics[width=0.8\textwidth]{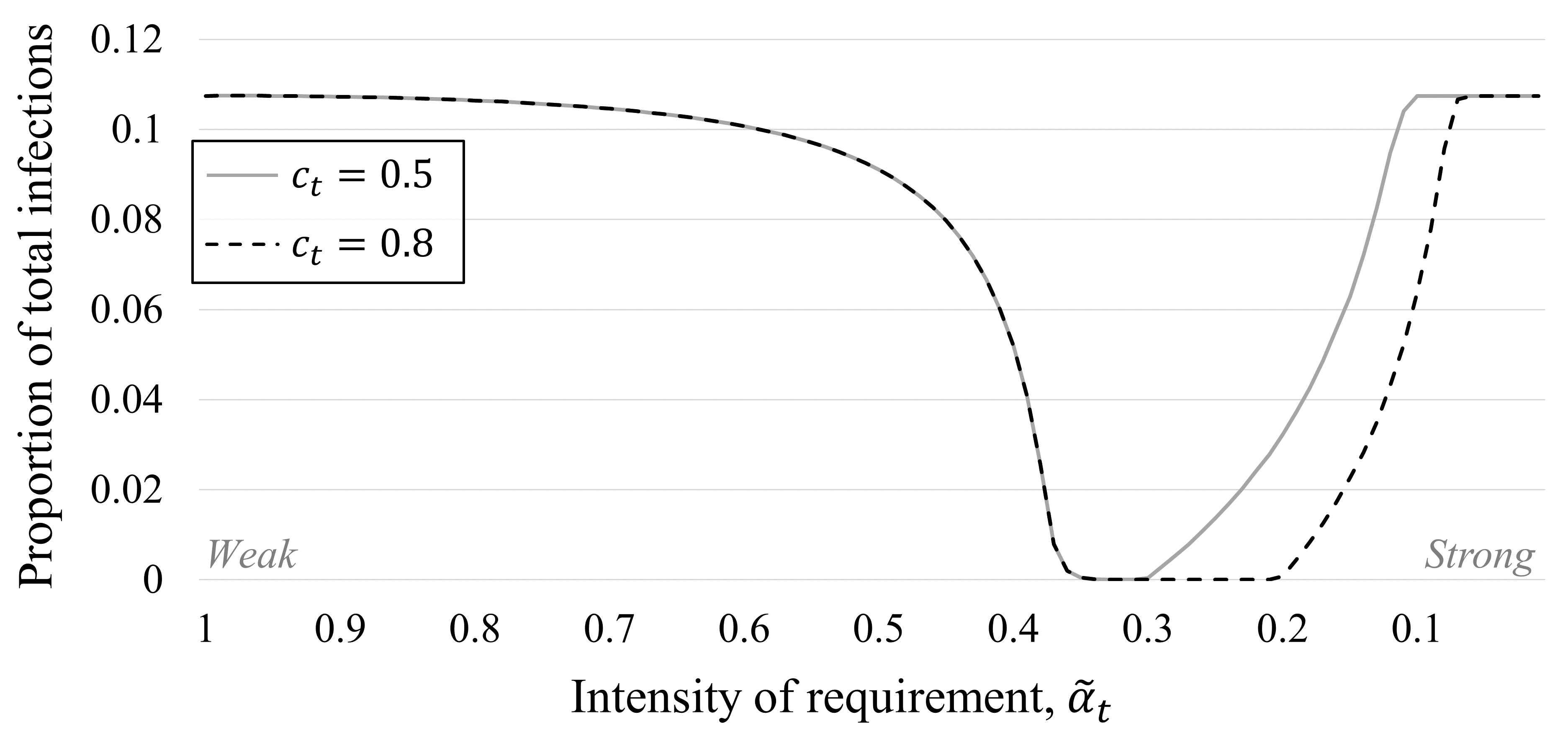}
\caption{Proportion of total infections for different noncompliance penalties: $c_t = 0.5$ and $0.8$.}
\label{fig:alpha-c-result}
{}
\end{figure}

\begin{figure}[htb!]
\centering
\includegraphics[width=0.75\textwidth]{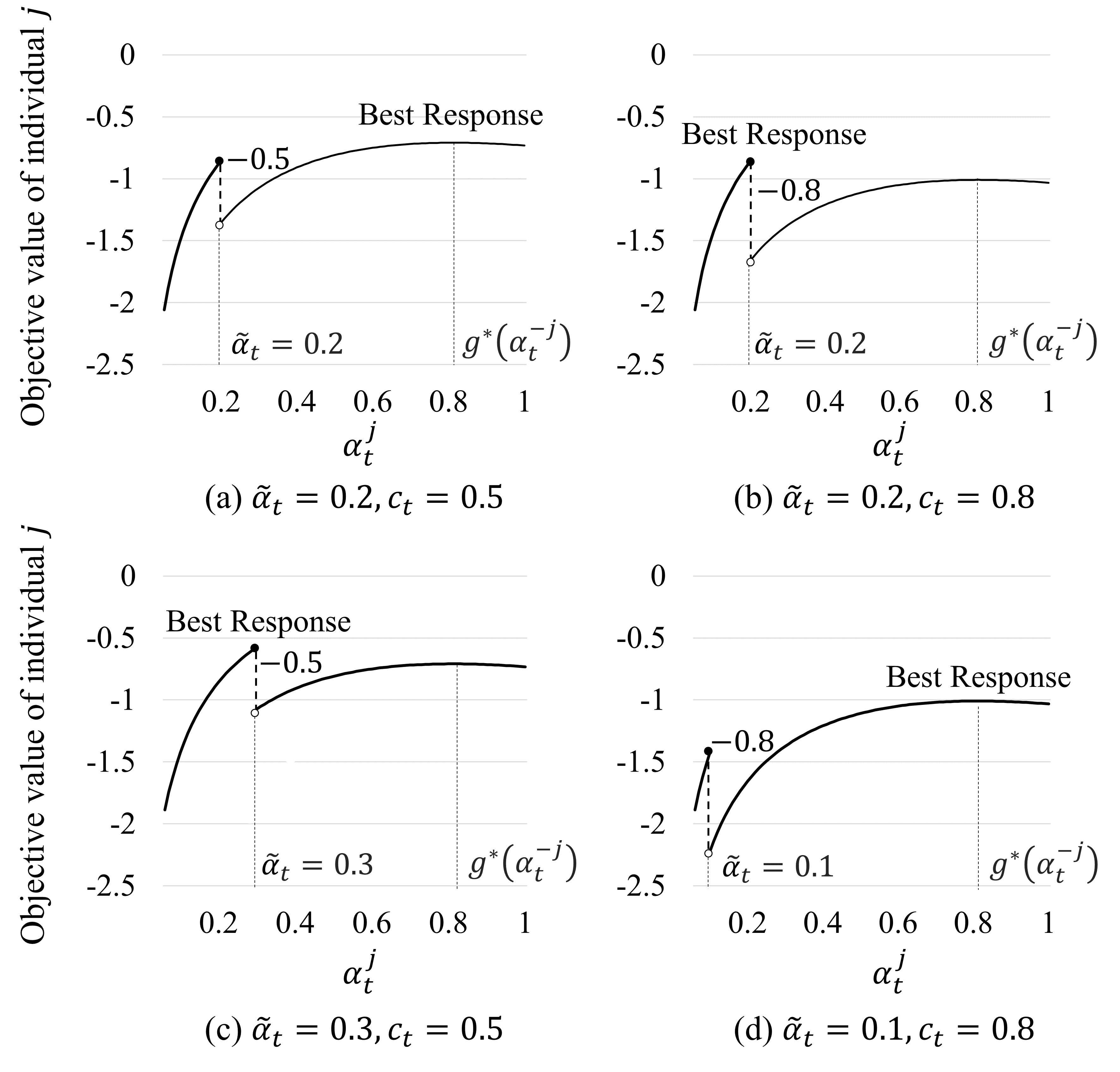}
\caption{Objective function as a function of $\alpha_t^j$ given a fixed $\alpha_t^{-j}$}
\label{fig:c-discussion}
{}
\end{figure}
\noindent Figure \ref{fig:alpha-c-result} shows the results when noncompliance penalty $c_t$ is 0.5 and 0.8 respectively. 
In this figure, we observe three regions: $\tilde{\alpha}_t$ is between 1.0 and $\sim0.3$  where the two curves are almost identical, between 0.3 and $\sim0.1$ in which the curve for $c_t = 0.8$ lies below the other, and the remaining interval where they become identical again. In other words, the effect of imposing higher $c_t$ is the shift of a starting point for the rebound in the total infections curve. This suggests that with a larger penalty, policy compliance is maintained for social distancing policies with a higher intensity. The enhanced compliance is due to the role of noncompliance penalty in determining the best response for an individual agent, which is illustrated in Figure \ref{fig:c-discussion}.

Figure \ref{fig:c-discussion}(a) and (b) are the objective function for agent $j$ as a function of its social activity level while keeping that of other individuals constant. The two curves are drawn for $\tilde{\alpha}_t = 0.2$, but with different values of $c_t$, 0.5 and 0.8 respectively. As shown in the figure, when $c_t= 0.8$, the drop of the objective function at $\tilde{\alpha}_t$ is so large that the value at its apex is lower than at $\tilde{\alpha}_t$, making $\tilde{\alpha}_t$ the best response for $j$. In contrast, for $c_t = 0.5$, the objective function value is higher at its apex, and the best response becomes $g^*(\alpha^{-j}_t)$, which means they choose to engage in a higher level of social activities than mandated. Yielding identical total infections in the first region can be explained similarly. When $\tilde{\alpha}_t$ is greater than $\sim$0.3, a drop of 0.5 at $\tilde{\alpha}_t$ is sufficient to make $\tilde{\alpha}_t$ the best response, as shown in Figure \ref{fig:c-discussion}(c). Likewise, for $\tilde{\alpha}_t \leq 0.1$, Figure \ref{fig:c-discussion}(d) demonstrates that $c_t$ of 0.8 is not large enough to shift the curve down to make $\tilde{\alpha}_t$ the best response.

\section{Discussion}
\label{sec:disc}

\subsection{Efficiency of a social distancing policy}
\label{subsec:disc-efficiency}


In the previous section, we demonstrated that overly strict social distancing policies might lead to noncompliance, undermining effectiveness, while penalties can help enforce compliance. Implementing extremely high penalties could theoretically ensure compliance and halt disease spread but may not be desirable due to significant social and economic costs. Therefore, a balance is necessary.

To address this trade-off, we evaluate the efficiency of social distancing policies by comparing the reduction in infections to the loss in social activities relative to no policy. We define the efficiency of a social distancing policy $(\tilde{\alpha}_t,c_t)$ as $E_{(\tilde{\alpha}_t,c_t)}$:
\begin{align}
\label{eq:efficiency}
    E_{(\tilde{\alpha}_t,c_t)} &= \frac{\textrm{relative reduction in total infections}, \hat{Z}_{(\tilde{\alpha}_t,c_t)}}{\textrm{relative loss of social activities}, \hat{\Theta}_{(\tilde{\alpha}_t,c_t)}} \\
    &= \frac{\{Z_{(1,0)}-Z_{(\tilde{\alpha}_t,c_t)}\} \, / \,{Z_{(1,0)}}} { 
    \sum_{t} (\Theta_{t|(1,0)} - \Theta_{t|(\tilde{\alpha}_t,c_t)}) \,/\, \sum_{t}\Theta_{t|(1,0)}} \nonumber
\end{align}
where $\Theta_{t|(\tilde{\alpha}_t,c_t)}$ denotes the social activity level at time $t$ averaged over the entire population, which is $\alpha^S(t) s(t) + \alpha^I(t) i(t) + \alpha^R(t)r(t)$, and $Z_{(\tilde{\alpha}_t,c_t)}$ denotes the total infections under a social distancing policy ($\tilde{\alpha}_t$, $c_t$).

\begin{figure}[htb]
\centering
\includegraphics[width=0.75\textwidth]{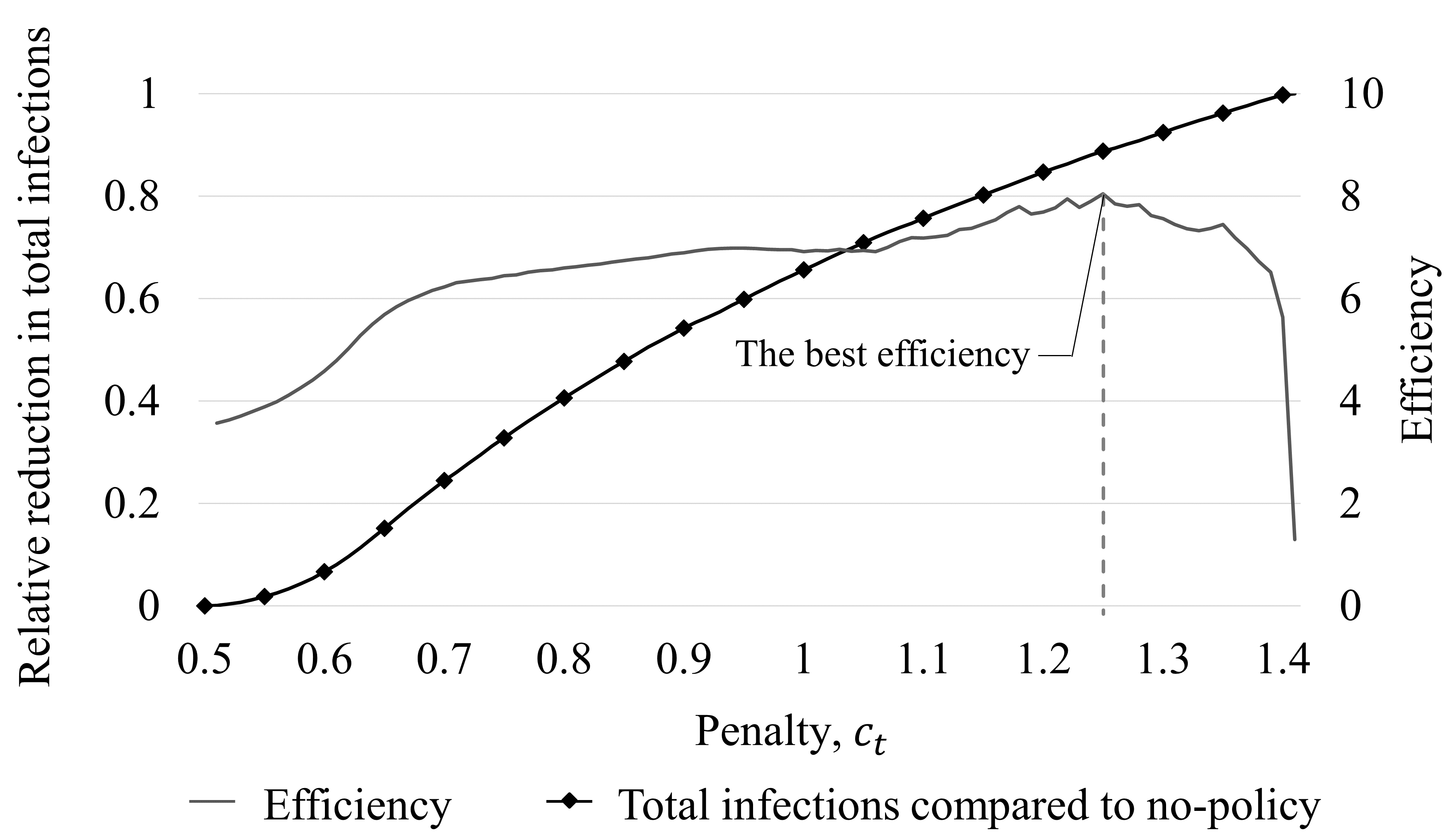}
\caption{Efficiency of social distancing policy and total infections as a function of penalty ($\tilde{\alpha}_t=0.1$)}
\label{fig:efficiency-result}
{}
\end{figure}

In Figure \ref{fig:efficiency-result}, we present results from implementing a stringent social distancing policy ($\tilde{\alpha}_t=0.1$) with varying penalties ($c_t$ from 0.5 to 1.4) to assess changes in efficiency $E_{(\tilde{\alpha}_t,c_t)}$ and reduction in infections $\hat{Z}_{(\tilde{\alpha}_t,c_t)}$. As expected, increasing $c_t$ leads to greater reductions in infections, achieving full reduction at $c_t = 1.4$. However, efficiency does not show a monotonic trend; it increases up to a point but begins to decline at $c_t = 1.25$ due to the faster growth in social activity loss compared to infection reduction. This indicates that the social and economic costs of social distancing might surpass its public health benefits, highlighting the need for careful evaluation of social distancing policies' efficiency.

\color{black}
\subsection{
Doing more with less: no social distancing for recovered population}
\label{disc:vaccine-pass}

From Figures \ref{fig:alpha-c-result} and \ref{fig:efficiency-result}, one might assume that higher penalties consistently reduce infections. However, the higher compliance from larger penalties does not always correlate with fewer infections.

\begin{figure}[htb]
\centering
\includegraphics[width=0.75\textwidth]{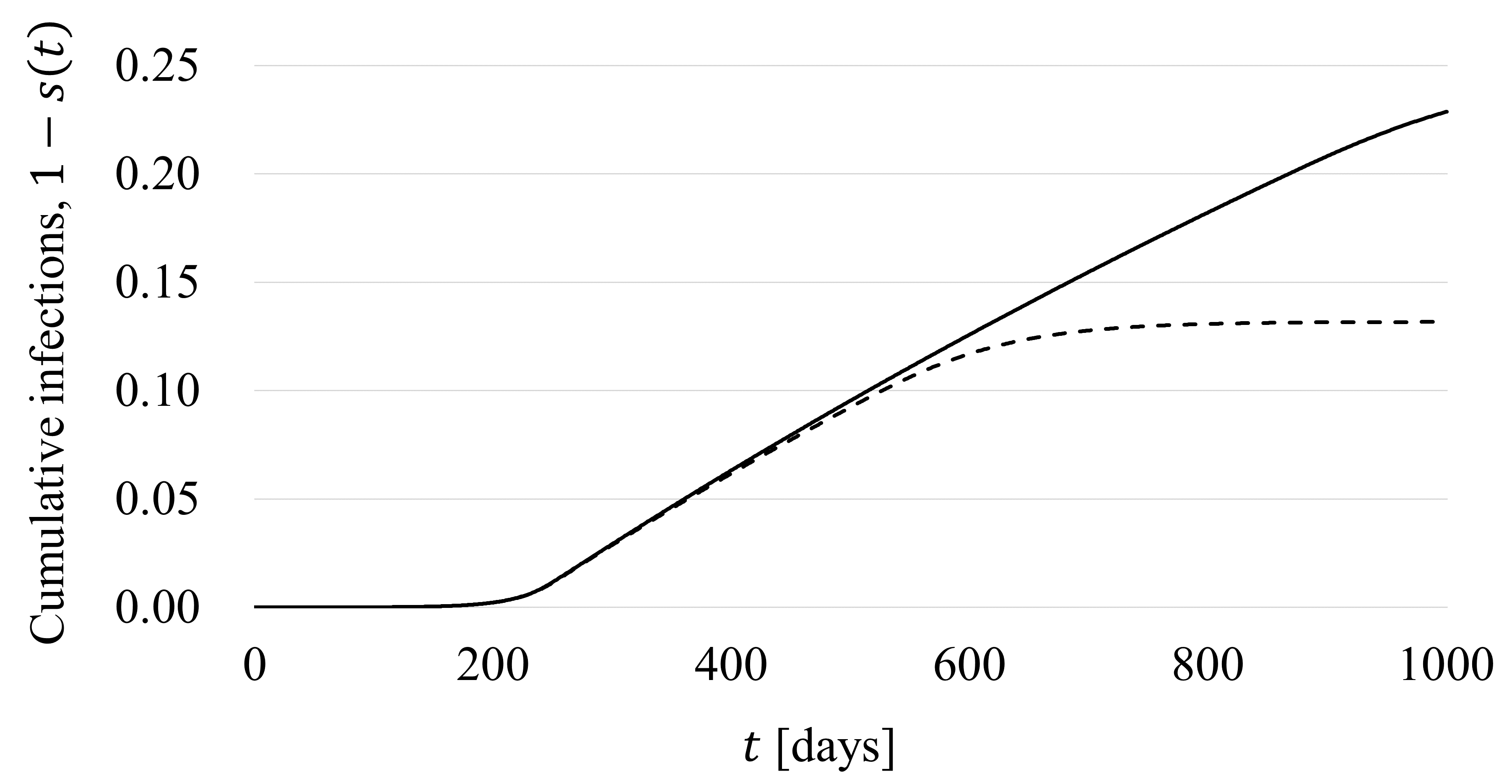}
\caption{Cumulative infections over time with $\tilde{\alpha}_t=0.4, c_t=0.4$ (solid); reduced penalty ($c_t=0.1$) applied to the recovered population (dotted)}
\label{fig:need-vaccine-pass}
\end{figure}

The solid line in Figure \ref{fig:need-vaccine-pass} depicts the cumulative infected population over time with a moderate penalty ($\tilde{\alpha}_t=0.4$ and $c_t=0.4$). For comparison, the dotted line shows outcomes with a reduced penalty ($c_t = 0.1$) applied specifically to the recovered population to encourage their noncompliance. Interestingly, infections are lower when the penalty for the recovered is decreased, an outcome that illustrates the complex dynamics of disease spread among different population groups.

This unexpected result emerges from the behavior of the recovered population. As seen in Eq.\eqref{eq:decision-R} and Figure \ref{graph:obj-R}, a lower penalty ($c_t$) leads the recovered to largely disregard the social distancing policy (\(\alpha^{*R}\) moves towards 1), increasing their social contacts. This uptick in social activity by the recovered dilutes their interactions with susceptible individuals, as indicated in Eq.\eqref{eq:adjusted_dsdt}, which paradoxically helps to control the spread of the disease. This suggests that exempting the recovered from strict social distancing mandates might enhance overall disease containment efforts.



\begin{figure}[htb]
\centering
\includegraphics[width=0.75\textwidth]{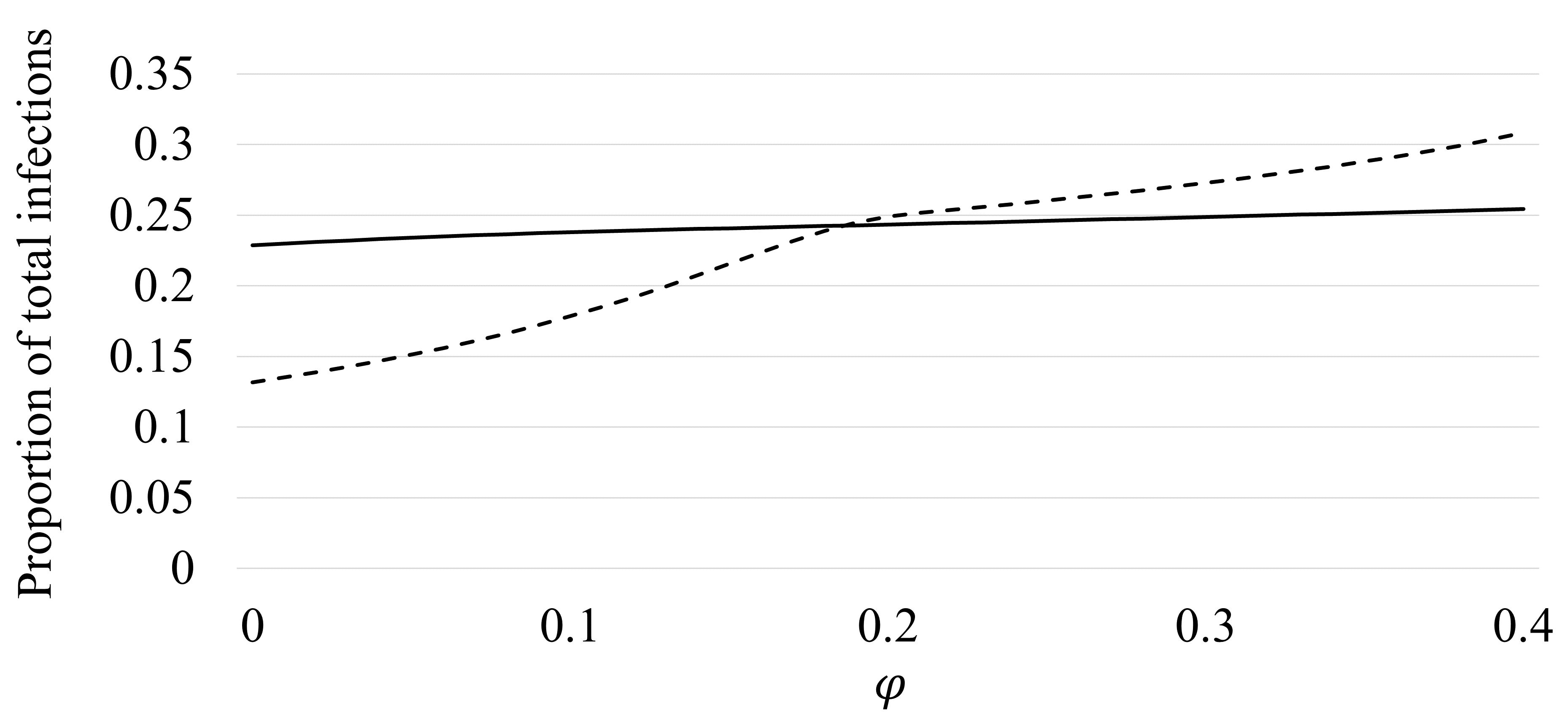}
\caption{Proportion of total infections as a function of reinfection rate $\varphi$:  uniform social distancing (solid), segmented social distancing (dotted)}
\label{fig:vaccine-pass}
{}
\end{figure}

The effectiveness of the above strategy hinges on the immunity of recovered individuals and the potential for reinfection. If reinfection risk is high, the strategy may fail. To assess this, we added a \textit{reinfected} compartment to our model, allowing for transitions from recovered to reinfected. The extended model details are in Appendix C. We test varying reinfection rates, $\varphi$, from 0.0 to 0.4, and present the findings in Figure \ref{fig:vaccine-pass}.

In Figure \ref{fig:vaccine-pass}, we analyze total infections under two scenarios: one with uniform social distancing across all population groups, and another exempting the recovered population from social distancing. For reinfection rates below 0.2, exempting the recovered reduces total infections. Above a 0.2 reinfection rate, uniform distancing proves more effective. The viability of segmented distancing hinges on the reinfection risk among the recovered. Given COVID-19's estimated reinfection rate below 0.2, as per \cite{pilz2022sars-reinfection}, a segmented approach could have been beneficial during the pandemic.

During COVID-19, regions like Hong Kong and South Korea implemented a \textit{vaccine pass} policy, granting privileges like access to supermarkets, restaurants, and theaters exclusively to vaccinated individuals, while the non-vaccinated faced strict distancing rules. The policy aimed to incentivize vaccination and stimulate economic activity. Our findings indicate that if vaccine efficacy is high, such policies can also curb disease spread, suggesting that segmented social distancing approaches merit consideration in policy design.


\section{Conclusion}
\label{sec:conclusion}

Public health authorities implement social distancing measures during disease outbreaks to mitigate spread. This research explores how adherence to such policies influences infectious disease proliferation, highlighting often overlooked noncompliance issues. We introduce a novel framework that incorporates individual compliance decisions within the dynamics of infectious diseases and policy design.

We model these decisions through a multi-agent ``social distancing game," determining a unique Nash equilibrium for each policy scenario. Using COVID-19 data from South Korea, numerical experiments evaluate compliance and its impact on disease transmission under various policies.

Our findings indicate that while social distancing generally helps control disease spread, overly strict measures can provoke noncompliance, reducing policy effectiveness. If severe measures are necessary, significant penalties may be required to ensure compliance. Moreover, the cost of social distancing, especially in severe forms like lockdowns, can outweigh benefits, necessitating careful consideration of policy efficiency. Additionally, for diseases with low reinfection risks, a segmented policy that exempts immune individuals could lower both infection rates and socioeconomic costs. \textcolor{black}{Most importantly, our study underscores the significance of considering potential noncompliance in social distancing policy design, a factor often ignored in the existing literature. We integrate a social distancing game with government-imposed policies to model noncompliance, highlighting that noncompliance penalties are crucial yet frequently overlooked policy elements.}

Before closing, we suggest potential areas for future research. First, the epidemic model in this study could be enhanced to better reflect the complexities of the real environment, such as integrating vaccination campaigns into the dynamics. Another improvement involves introducing population heterogeneity, considering individual differences like age, gender, and health conditions to refine the analysis of social distancing policies. \textcolor{black}{Specifically, one may incorporate age heterogeneity by adopting age group specific transmission rate, cost of infection, sensitivity to social activities.  Incorporating heterogeneity certainly brings computational challenges, requiring development of algorithmic solution approached. The Nash equilibrium discussed in this paper would then serve as an effective means to assess the algorithm's performance. Additionally, refining the utility function to reflect more complex decision-making, where the socialization utility depends on the activity levels of others, would increase model relevance.} Lastly, expanding this work to design dynamic social distancing policies as an optimization problem over the disease's course could offer more strategic insights, moving beyond assessing static policies to find optimal time-varying strategies.

\color{black}
\section*{Acknowledgment}
This work was supported by the National Research Foundation of Korea(NRF) grant funded by the Korea government (MSIT) (No. NRF-2022R1A2B5B01002615).

\section*{CRediT authorship contribution statement}
\noindent \textbf{Hyelim Shin}: Conceptualization, Methodology, Software, Investigation, Visualization, Writing - Original Draft. \textbf{Taesik Lee}: Conceptualization, Writing - Review \& Editing, Supervision.

\section*{Declaration of Generative AI and AI-assisted technologies in the writing process}
During the preparation of this work the author(s) used ChatGPT in order to refine some part of the text. After using this tool/service, the author(s) reviewed and edited the content as needed and take(s) full responsibility for the content of the publication.

\newpage
\setcounter{page}{1}

\centerline{\huge Supplementary Materials}

\appendix

\section{Optimality and Uniqueness of a Nash equilibrium}
\label{A:solution-uniqueNE}

In Section \ref{subsec:solution-SI}, we classify the scenarios into four distinct cases based on the given policy parameters ($\tilde{\alpha}$ and $c$) and the epidemic state variables $(\hat{s},\hat{i}, \hat{q},$ and $\hat{r})$. For each case, we compute a Nash equilibrium. In this section, we demonstrate that these cases are mutually exclusive, and each case has a unique Nash equilibrium. Proposition \ref{thm:BR} guarantees the optimality of each Nash equilibrium and Proposition \ref{thm:hat-btw-tilde&bar} shows the existence of the mixed Nash equilibrium, $\hat{\alpha}$, in case 4. Proposition \ref{thm:if-case2-not-case3} and \ref{thm:if-case3-not-case2} shows that there is no intersection between any two cases. And, we show the uniqueness of Nash equilibrium for each case through Proposition \ref{thm:uniqueNE-case1}-\ref{thm:uniqueNE-case4}.

\color{black}
\begin{repproposition}{thm:BR}

The best response of individual $j$ for a given strategy of other players is
\begin{equation*}
    BR(\alpha^{-j}) = \begin{cases}
                        \tilde{\alpha} \text{ with probability } p\\
                        g^*( \alpha^{-j}) \text{ with probability } 1-p\\
                      \end{cases}
\end{equation*}
for some $p \in [0,1]$. 
\end{repproposition}

\begin{proof}
\begin{align*}
    \text{Let } &f_1(\alpha^j):=\log{\alpha^j}-\alpha^j+1,\text{ and}\\
    &f_2(\alpha^j,\alpha^{-j}):=\frac{\alpha^j i \alpha^{-j} \pi \kappa d}{s\alpha^{-j}+i\alpha^{-j}+r\alpha^{*R}}.
\end{align*}
Then,
\begin{equation*}
    u(\alpha^j,\alpha^{-j}) = \begin{cases}
        f_1(\alpha^j)+f_2(\alpha^j,\alpha^{-j}) & \alpha^j \leq \tilde{\alpha}\\
        f_1(\alpha^j)+f_2(\alpha^j,\alpha^{-j}) -c & \alpha^j > \tilde{\alpha}\\
    \end{cases}
\end{equation*}
By the definition of $g^*( \alpha^{-j})$, $f_1(g^*( \alpha^{-j}))+f_2(g^*( \alpha^{-j}),\alpha^{-j}) \geq f_1(\alpha^j)+f_2(\alpha^j,\alpha^{-j}) \ \forall \alpha^j \in [0,1]$. And $f_1(g^*( \alpha^{-j}))+f_2(g^*( \alpha^{-j}),\alpha^{-j}) = f_1(\alpha^j)+f_2(\alpha^j,\alpha^{-j})$ if and only if $\alpha^j=g^*( \alpha^{-j})$.

\noindent i) $\tilde{\alpha} > g^*( \alpha^{-j})$

For $\alpha^j \in [0,\tilde{\alpha}]$,
\begin{align*}
u(g^*( \alpha^{-j}), \alpha^{-j}) &= f_1(g^*( \alpha^{-j}))+f_2(g^*( \alpha^{-j}),\alpha^{-j})\\
&\geq f_1(\alpha^j)+f_2(\alpha^j,\alpha^{-j})\\
&= u(\alpha^j, \alpha^{-j}).
\end{align*}

For $\alpha^j \in (\tilde{\alpha},1]$,
\begin{align*}
u(g^*( \alpha^{-j}), \alpha^{-j}) &=f_1(g^*( \alpha^{-j}))+f_2(g^*( \alpha^{-j}),\alpha^{-j}) \\
&> f_1(\alpha^j)+f_2(\alpha^j,\alpha^{-j}) \\
&> f_1(\alpha^j)+f_2(\alpha^j,\alpha^{-j}) -c \\
&= u(\alpha^j, \alpha^{-j}).
\end{align*}

As $u(g^*( \alpha^{-j}), \alpha^{-j})$ is the maximum value, $g^*( \alpha^{-j})$ is the only best response of individual $j$. In this case, $p=0$.

\noindent ii) $\tilde{\alpha} < g^*( \alpha^{-j})$

For $\alpha^j \in [0, \tilde{\alpha}]$, $f_1(\alpha^j)+f_2(\alpha^j,\alpha^{-j})$ is monotonically increasing in $[0, \tilde{\alpha}]$. Therefore, $u(\alpha^j, \alpha^{-j})$ attains its maximum value at $\alpha^j=\tilde{\alpha}$ within the interval $[0, \tilde{\alpha}]$.

For $\alpha^j \in (\tilde{\alpha},1]$,
\begin{align*}
u(g^*( \alpha^{-j}), \alpha^{-j}) &=f_1(g^*( \alpha^{-j}))+f_2(g^*( \alpha^{-j}),\alpha^{-j})-c \\
&\geq f_1(\alpha^j)+f_2(\alpha^j,\alpha^{-j}) -c \\
&= u(\alpha^j,\alpha^{-j})
\end{align*}
by definition of $g^*( \alpha^{-j})$. Then, the best response is as follows for action $\{\tilde{\alpha}, g^*( \alpha^{-j})\}$.

\begin{equation*}
    BR(\alpha^{-j}) = \begin{cases}
                        (1,0) &\text{if } u(\tilde{\alpha}, \alpha^{-j}) > u(g^*( \alpha^{-j}), \alpha^{-j})\\
                        (p,1-p):0\leq p \leq 1 &\text{if } u(\tilde{\alpha}, \alpha^{-j}) = u(g^*( \alpha^{-j}), \alpha^{-j}) \\
                        (0,1) &\text{if } u(\tilde{\alpha}, \alpha^{-j}) < u(g^*( \alpha^{-j}), \alpha^{-j})\\
                      \end{cases}
\end{equation*}

\end{proof}

\textcolor{black}{
\begin{proposition}
\label{thm:hat-btw-tilde&bar}
Suppose $\tilde{\alpha}<\bar{\alpha}$. Then, if $u(\bar{\alpha},\bar{\alpha})<u(\tilde{\alpha},\bar{\alpha})$ and $u(\tilde{\alpha},\tilde{\alpha})<u(g^*(\tilde{\alpha}),\tilde{\alpha})$, there exists $\alpha \in (\tilde{\alpha},\bar{\alpha})$ that satisfies $u(\tilde{\alpha},\alpha)=u(g^*(\alpha),\alpha)$ and $\alpha < g^*(\alpha)$.
\end{proposition}
}

\begin{proof}

\textcolor{black}{
Let $U(\alpha):=u(\tilde{\alpha},\alpha)-u(g^*(\alpha),\alpha)$. We first prove that there exists $\alpha \in (\tilde{\alpha},\bar{\alpha})$ such that $U(\alpha)=0$ by the intermediate value theorem. In other words, if (i) $U(\alpha)$ is continuous on $[\tilde{\alpha},\bar{\alpha}]$ and (ii) $U(\tilde{\alpha})U(\bar{\alpha})<0$, the existence of such $\alpha$ is guaranteed.}

\textcolor{black}{
(i) \textit{Continuity}
By Eq.\eqref{eq:utilityfcn}, $u(\alpha^j,\alpha^{-j})$ is continuous in $\alpha^{-j}$ and $u(\alpha^j,\alpha^{-j})$ is continuous on $(\tilde{\alpha},1]\times [\tilde{\alpha},\bar{\alpha}]$. Trivially the first term of $U(\alpha)$, $u(\tilde{\alpha},\alpha)$, is continuous. For the continuity of the second term, it is enough to show that $g^*(\alpha)>\tilde{\alpha}$ for all $\alpha \in [\tilde{\alpha},\bar{\alpha}]$ since $g^*(\alpha)$ is continuous in $\alpha$ by Eq.\eqref{eq:g*}. The condition $\tilde{\alpha}\leq \alpha \leq \bar{\alpha}$ guarantees $g^*(\bar{\alpha}) \leq g^*(\alpha) \leq g^*(\tilde{\alpha})$ because $g^*$ is a decreasing function, and $g^*(\bar{\alpha})=\bar{\alpha}$. Combining the two inequalities, $\tilde{\alpha}<g^*(\alpha)$ by the given assumption $\tilde{\alpha}<\bar{\alpha}$. Therefore, the function $U(\alpha)$ is continuous on $[\tilde{\alpha},\bar{\alpha}]$.}

\textcolor{black}{
(ii) \textit{Negative product}
By the assumption -- $u(\bar{\alpha},\bar{\alpha})<u(\tilde{\alpha},\bar{\alpha})$ and $u(\tilde{\alpha},\tilde{\alpha})<u(g^*(\tilde{\alpha}),\tilde{\alpha})$ --, $U(\tilde{\alpha})$ is negative and $U(\bar{\alpha})$ is positive, yielding that  their product $\ U(\tilde{\alpha})U(\bar{\alpha})$ is negative.}

\textcolor{black}{In addition, $\alpha < \bar{\alpha}=g^*(\bar{\alpha})<g^*(\alpha)$ since $g^*$ is a decreasing function.}
\end{proof}

\color{black}

According to the definition of the cases, there is no intersection between any two cases, except for case 2 and case 3. Therefore, if case 2 and case 3 do not coincide with each other, it can be concluded that the four cases are entirely separate and do not intersect.

\begin{proposition}
\label{thm:if-case2-not-case3}
For given $\tilde{\alpha}$ and $c$, if an epidemic state $(\hat{s},\hat{i}, \hat{q},$ and $\hat{r})$ is included in case 2, then the state is not included in case 3.
\end{proposition}
\begin{proof}
What should be proven is that

\begin{equation*}
    \text{where }  \tilde{\alpha}<\bar{\alpha}, \ u(\bar{\alpha},\bar{\alpha}) \geq u(\tilde{\alpha},\bar{\alpha}) \implies u(\tilde{\alpha},\tilde{\alpha}) < u(g^*(\tilde{\alpha}),\tilde{\alpha}) .
\end{equation*}

If $u(\tilde{\alpha},\tilde{\alpha}) < u(\bar{\alpha},\tilde{\alpha})$, $u(\tilde{\alpha},\tilde{\alpha}) < u(g^*(\tilde{\alpha}),\tilde{\alpha})$ trivially by Definition \ref{def:g*}. Let $f(\alpha):=\log(\alpha)-\alpha+1$.

\begin{gather*}
    u(\bar{\alpha},\bar{\alpha}) \geq u(\tilde{\alpha},\bar{\alpha}) \Leftrightarrow f(\bar{\alpha})-f(\tilde{\alpha})-c \geq \kappa\pi d \frac{\bar{\alpha}\hat{i}}{\bar{\alpha}\hat{s}+\bar{\alpha}\hat{i}+\alpha^{*R}\hat{r}} (\bar{\alpha}-\tilde{\alpha})\\
    \text{and}\\
    u(\bar{\alpha},\tilde{\alpha}) > u(\tilde{\alpha},\tilde{\alpha}) \Leftrightarrow f(\bar{\alpha})-f(\tilde{\alpha})-c > \kappa\pi d \frac{\tilde{\alpha}\hat{i}}{\tilde{\alpha}\hat{s}+\tilde{\alpha}\hat{i}+\alpha^{*R}\hat{r}} (\bar{\alpha}-\tilde{\alpha}).\\
    \therefore \kappa\pi d \frac{\bar{\alpha}\hat{i}}{\bar{\alpha}\hat{s}+\bar{\alpha}\hat{i}+\alpha^{*R}\hat{r}} (\bar{\alpha}-\tilde{\alpha}) > \kappa\tau d \frac{\tilde{\alpha}\hat{i}}{\tilde{\alpha}\hat{s}+\tilde{\alpha}\hat{i}+\alpha^{*R}\hat{r}} (\bar{\alpha}-\tilde{\alpha}) \implies u(\bar{\alpha},\tilde{\alpha}) > u(\tilde{\alpha},\tilde{\alpha}).
\end{gather*}

\begin{align*}
    \tilde{\alpha}<\bar{\alpha} &\Leftrightarrow \frac{\bar{\alpha}\hat{i}}{\bar{\alpha}\hat{s}+\bar{\alpha}\hat{i}+\alpha^{*R}\hat{r}}  >  \frac{\tilde{\alpha}\hat{i}}{\tilde{\alpha}\hat{s}+\tilde{\alpha}\hat{i}+\alpha^{*R}\hat{r}}\\
     &\Leftrightarrow \kappa\pi d \frac{\bar{\alpha}\hat{i}}{\bar{\alpha}\hat{s}+\bar{\alpha}\hat{i}+\alpha^{*R}\hat{r}} (\bar{\alpha}-\tilde{\alpha}) > \kappa\pi d \frac{\tilde{\alpha}\hat{i}}{\tilde{\alpha}\hat{s}+\tilde{\alpha}\hat{i}+\alpha^{*R}\hat{r}} (\bar{\alpha}-\tilde{\alpha}).
\end{align*}
\end{proof}

\begin{proposition}
\label{thm:if-case3-not-case2}
For given $\tilde{\alpha}$ and $c$, if epidemic state $(\hat{s},\hat{i},\hat{r})$ is included in case 3, then the state is not included in case 2.
\end{proposition}
\begin{proof}
It is simply proved similar to proof of Proposition \ref{thm:if-case2-not-case3}.
\end{proof}

\begin{proposition}
\label{thm:uniqueNE-case1}
$\bar{\alpha}$ is the unique Nash equilibrium under case 1.
\end{proposition}

\begin{proof}
i) $\tilde{\alpha} (\neq \bar{\alpha})$ is not a Nash equilibrium.
\begin{align*}
    \tilde{\alpha} > \bar{\alpha} &\Leftrightarrow g^*(\bar{\alpha}) > g^*(\tilde{\alpha})  &&\because \ g^*(\alpha) \text{ is decreasing function of } \alpha \\
    &\Leftrightarrow \tilde{\alpha} > \bar{\alpha} > g^*(\tilde{\alpha}) &&\because \ g^*(\bar{\alpha})=\bar{\alpha}\\
    \therefore u(g^*&(\tilde{\alpha}),\tilde{\alpha})>u(\tilde{\alpha},\tilde{\alpha}).
\end{align*}
$\tilde{\alpha} (\neq \bar{\alpha})$ is not a best response where $\alpha^{-j}=\tilde{\alpha}$.

ii) Mixed strategy $\hat{\alpha}$ is not a Nash equilibrium.

Let assume that there exists a mixed strategy $\hat{\alpha}$ that satisfies \eqref{eq:mu-hat-def} and \eqref{eq:mu-hat-cond}.

\begin{align*}
    \tilde{\alpha} \geq \bar{\alpha} &\Leftrightarrow g^*(\bar{\alpha}) \geq g^*(\tilde{\alpha})  &&\because \ g^*(\alpha) \text{ is decreasing function of } \alpha \\
    &\Leftrightarrow \tilde{\alpha} \geq \bar{\alpha} \geq g^*(\tilde{\alpha}) &&\because \ g^*(\bar{\alpha})=\bar{\alpha}\\
    &\Leftrightarrow g^*(\hat{\alpha}) > \tilde{\alpha} \geq \bar{\alpha} \geq g^*(\tilde{\alpha})  &&\because \ g^*(\hat{\alpha}) > \tilde{\alpha} \text{ in } \eqref{eq:mu-hat-def}.
\end{align*}
However, $g^*(\hat{\alpha})>g^*(\tilde{\alpha}) \Leftrightarrow \tilde{\alpha}>\hat{\alpha}$. It is contradiction to \eqref{eq:mu-hat-def}.
\end{proof}

\begin{proposition}
\label{thm:uniqueNE-case2}
$\bar{\alpha}$ is the unique Nash equilibrium under case 2.
\end{proposition}
\begin{proof}
By Proposition \ref{thm:if-case2-not-case3}, $\tilde{\alpha}$ is not a Nash equilibrium. We need to prove that $\hat{\alpha}$ is not a Nash equilibrium under case 2.

Let assume that there exists a mixed strategy $\hat{\alpha}$ such that satisfies \eqref{eq:mu-hat-def} and \eqref{eq:mu-hat-cond}.

\begin{align*}
    \text{Let } &U(\alpha):=u(\tilde{\alpha},\alpha)-u(g^*(\alpha),\alpha), \text{ for $\alpha$ s.t. } g^*(\alpha)>\tilde{\alpha}\\ &f(\alpha):=\log{\alpha}-\alpha+1,\text{ and}\\
    &g(\alpha):=\frac{\alpha i \pi \kappa d}{s\alpha+i\alpha+r\alpha^{*R}}.
\end{align*}

\begin{align*}
    U(\alpha) &= f(\tilde{\alpha})-\tilde{\alpha}g(\alpha)-f(g^*(\alpha))+g^*(\alpha)g(\alpha)+c.\\
    \frac{d}{d\alpha}U(\alpha) &= -\tilde{\alpha}g'(\alpha)-f'(g^*(\alpha))(\frac{d}{d\alpha}g^*(\alpha))+(\frac{d}{d\alpha}g^*(\alpha))g(\alpha)+g^*(\alpha)g'(\alpha)\\
    &=g'(\alpha)(g^*(\alpha)-\tilde{\alpha})+(\frac{d}{d\alpha}g^*(\alpha))(g(\alpha)-f'(g^*(\alpha)))\\
    &=g'(\alpha)(g^*(\alpha)-\tilde{\alpha}) \ \ \because \  g(\alpha)=f'(g^*(\alpha)) \text{ by Definition } \ref{def:g*}.
\end{align*}
As $g'(\alpha)>0$, $U(\alpha)$ is nondecreasing function for $\alpha$ such that $g^*(\alpha)>\tilde{\alpha}$. In case 2, $\tilde{\alpha}<\hat{\alpha}<\bar{\alpha}=g^*(\bar{\alpha})<g^*(\hat{\alpha})<g^*(\tilde{\alpha})$ by \eqref{eq:mu-hat-def}. Since $U(\tilde{\alpha})<0$ (by Proposition \ref{thm:if-case2-not-case3}) and $U(\bar{\alpha}) \leq 0$ (by condition of case 2), $U(\hat{\alpha})<0$. It is a contradiction to \eqref{eq:mu-hat-cond}.
\end{proof}

\begin{proposition}
\label{thm:uniqueNE-case3}
$\tilde{\alpha}$ is the unique Nash equilibrium under case 3.
\end{proposition}
\begin{proof}
It is simply proved similar to proof of Proposition \ref{thm:uniqueNE-case2}.
\end{proof}

\begin{proposition}
\label{thm:uniqueNE-case4}
$\hat{\alpha}$ is the unique Nash equilibrium under case 4.
\end{proposition}
\begin{proof}
There does not exists $\hat{\alpha}$ such that $g^*(\hat{\alpha}) \leq \tilde{\alpha}$ because it is contradiction to \eqref{eq:mu-hat-def}.

In the proof of Proposition \ref{thm:uniqueNE-case2}, $U(\alpha)$ is nondecreasing function for $\alpha$ such that $g^*(\alpha)>\tilde{\alpha}$. In case 4, $\tilde{\alpha}<\hat{\alpha}<\bar{\alpha}=g^*(\bar{\alpha})<g^*(\hat{\alpha})<g^*(\tilde{\alpha})$ by \eqref{eq:mu-hat-def}. $\hat{\alpha}$ is the unique solution of $U(\alpha)=0$.
\end{proof}

\color{black}
\section{Additional experimental results}
\label{B:results}

\begin{figure}[htb!]
\centering
\includegraphics[width=0.8\textwidth]{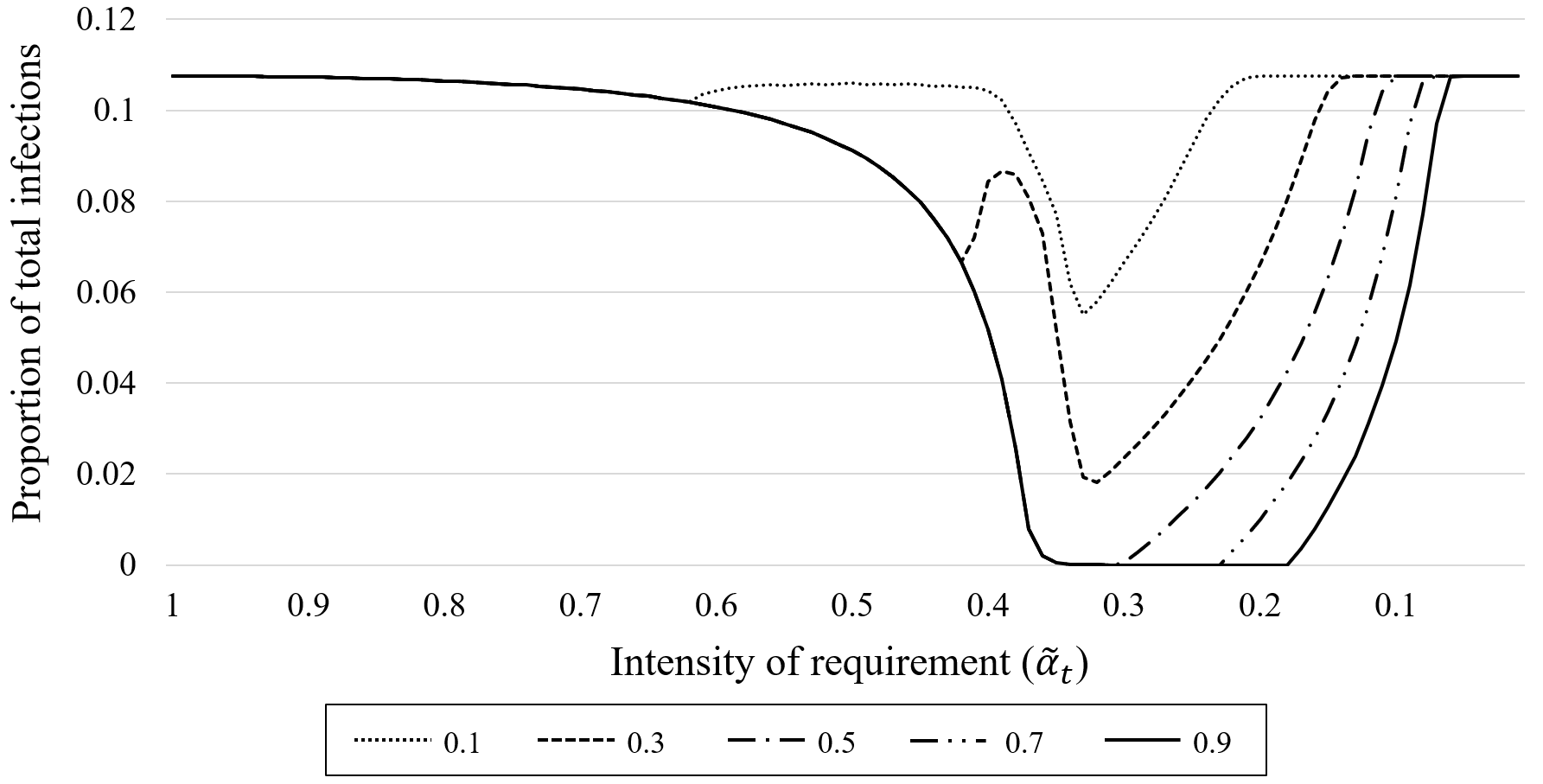}
\caption{Proportion of the total infected population according to social distancing intensity for various values of penalty ($c_t$).}
\label{fig:alpha-tilde-result-appendix}
{}
\end{figure}

Figure \ref{fig:alpha-tilde-result-appendix} depicts total infections according to social distancing intensity $\tilde{\alpha}_t$ for various penalty value $c_t$. The figure illustrates the prevalent occurrence of rebounds in total infections. Furthermore, it demonstrates that increasing the value of $c_t$ leads to a shift in the rebound's starting point and a decrease in the lowest point, indicating the impact of higher penalty values on the infection dynamics.

\section{SIQR model with reinfection}
\label{C:reinfection}

\begin{figure}[htb]
\centering
\includegraphics[width=\textwidth]{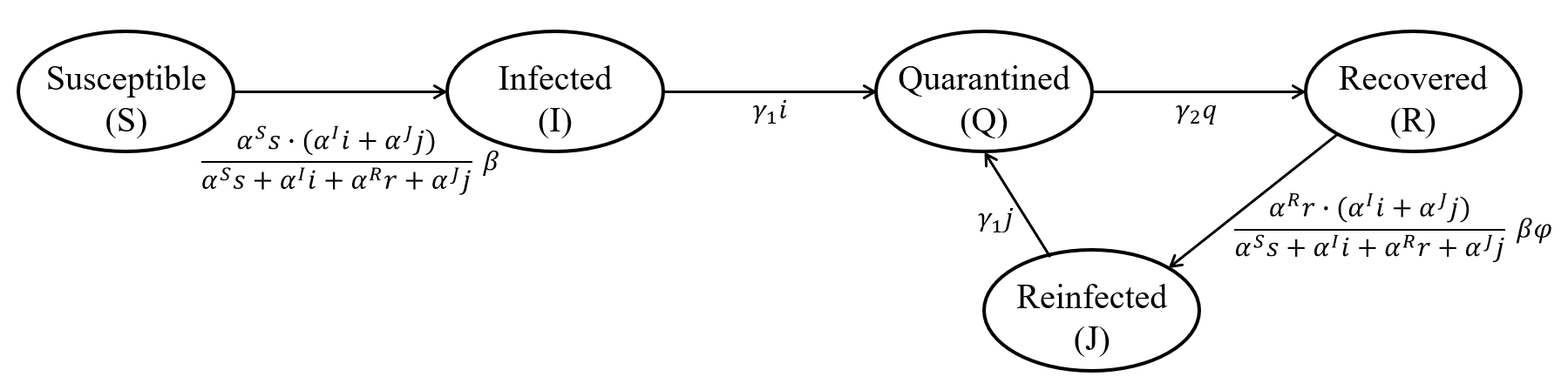}
\caption{Illustration of an epidemic model with reinfection rate}
\label{fig:SIRJ}
{}
\end{figure}

We add reinfected population to the proposed model in section \ref{sec:problem} to express reinfection. There are three differences between reinfection and infection. First, $\varphi$ makes the difference between the infection rate and the reinfection rate. $\varphi$ denotes the transmission probability given contact with a recovered individual and an infected individual compared to contact with a susceptible individual and an infected individual. If $\varphi$ is 0, the new epidemic model is the same as the existing model because reinfection does not occur at all. In general, the recovered population has a lower transmission probability than the susceptible population, so $\varphi$ is assumed to be less than 1. Second, the reinfected population has the same social activity level as the recovered population because the reinfected population does not realize that they are infected again, just as the infected population has the same social activity level as the susceptible population. That is, the reinfected population has a higher activity level than the infected population. The last difference is the death cost. A recovered individual is less likely to lead to death even if reinfected, so death due to reinfection is not considered for the utility function.





\afterpage{
\bibliography{arxiv}
}

\end{document}